\documentclass[11pt,notitlepage,twoside,a4paper]{article}
\usepackage{graphicx}
\usepackage{psfrag}
\usepackage{amssymb}
\usepackage{amsmath,amsthm}
\usepackage{graphicx}
\usepackage{amsmath, amssymb, amsfonts,enumerate}
\usepackage[french]{babel}
\usepackage{t1enc}
\usepackage{psfrag}
\usepackage{mathrsfs}
\newcommand{\summ}[3]{\setlength{\unitlength}{1cm}
\raisebox{0.0cm}[0.65cm][0.5cm]{
\begin{picture}(0.75,0.65)\put(0,0){\makebox{$\displaystyle\sum\!^*$}}
\put(0.25,-0.36){\makebox(0,0){$\scriptstyle #1=#2$}}
\put(0.25,0.53){\makebox(0,0){$\scriptstyle #3$}}
\end{picture}}   }
\newcommand {\hl }{\\\hline}
\newcommand {\ba }{\begin{array}}
\newcommand {\ea }{\end{array}}

\newtheorem{defi}{Definition}[section]
\newtheorem{theo}[defi]{Theorem}
\newtheorem{propo}[defi]{Proposition}
\newtheorem{lem}[defi]{Lemma}

\author{Hocine Sellama\\ \small IRMA - UMR 7501 CNRS/ULP\\\small 7 rue René
  Descartes - 67084 Strasbourg Cedex, France\\\small email: 
sellama@math.u-strasbg.fr}
\title{\bf On the distance between separatrices for the  discretized logistic differential equation }
\date{}
\numberwithin{equation}{section}
\begin{document}
\maketitle
 \renewcommand{\abstractname}{abstract}
 \begin{abstract}
 In this paper, we consider the discretization
\begin{equation*}
y(t+\varepsilon)=y(t-\varepsilon)+2\varepsilon\big(1-y(t)^{2}\big),
\end{equation*}
$\varepsilon>0$ a small parameter, of the logistic differential
 equation $y'=1-y^{2}$, which can  also be seen as
 a discretization of the system  
\begin{eqnarray*}
\begin{cases}
y'=2\big(1-v^{2}\big),\\
v'= 2\big(1-y^{2}\big).
\end{cases}
\end{eqnarray*}
This system has two saddle points at  $A=(1,1)$, $B=(-1, -1)$ and
there exist
 stable and unstable manifolds. We will show that the stable
manifold  $W_{s}^{+}$\, of the point $A=(1,1)$ and the unstable
manifold  $W_{i}^{-}$\, of the point $B=(-1, -1)$ for the discretization do not coincide. The vertical distance
between these two manifolds is exponentially small but not zero, in
particular we give an asymptotic estimate of this distance.  For
this purpose we will use a method adapted from the  paper of
  Sch\"afke-Volkmer \cite{SV} using formal series  and  accurate  estimates  of the coefficients.

 \end{abstract}


\section{Introduction}
\noindent We consider the logistic equation 
\begin{equation}
 y'=1-y(t)^{2},
\label{1.1}\end{equation} 
 whose solutions are $ y(t)=\tanh(t+c)$.
The discretization of this equation by 
Nystr\^om's method,  which consists in replacing the derivative by the
symmetrical difference, gives the recurrence:

\begin{equation}
 y_{n+1}= y_{n-1}+2\varepsilon\big(1-y_{n}^{2}\big).
\label{1.2}\end{equation}
 With the initial conditions \quad $y_{0}=0,\  y_{1}=\varepsilon,$ we calculate  the discrete solution  $y_{n}$ with 1600 iterates and
  $\varepsilon=\frac{1}{20}$ (see FIG.1).
\begin{figure}[!th]
   \centering
\includegraphics[height=9cm,width=9.5cm]{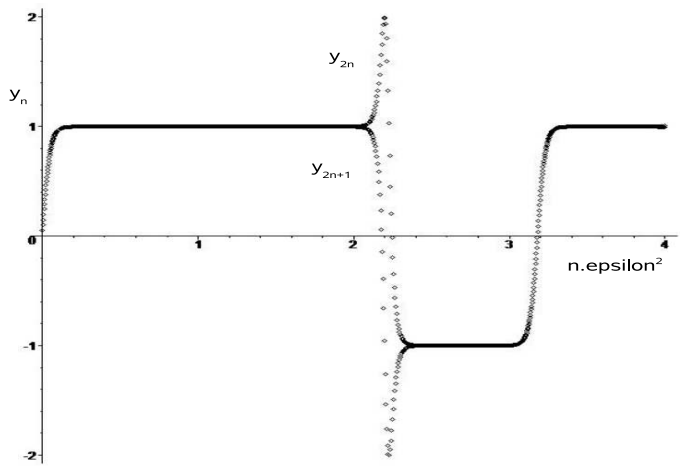}
  \caption{\small{This figure represents $y_{n}$  as a function of $n\varepsilon^{2} $.}}
  \label{p213}
\end{figure}
We observe that the discrete solution  rejoins quickly the level
$y=1$,  it stays close to this point  for a relatively long time  before it leaves its neighborhood, then for the
little while the even numbered points follow one curve and the odd
numbered follow another curve, the two curves meet at the level
$y=-1$, the discrete solution stays close to this point  for a 
relatively long time,  then   it remakes the same cycle.

Now letting $u_{n}=y_{2n},\,\,\, v_{n}=y_{2n+1} $, we obtain a recurrence of first order in the plane 
\begin{equation}
(u_{n+1},v_{n+1})=\Phi(u_{n},v_{n})
\label{1.3}\end{equation}
where the diffeomorphism  $\Phi:\mathrm{IR}^{2}\mapsto\mathrm{IR}^{2} $
is defined by 
\begin{eqnarray*}
u_{1}&=& u+2\varepsilon\big(1-v^{2}\big),\nonumber\\
v_{1}&=& v+2\varepsilon\big(1-u_{1}^{2}\big).
\end{eqnarray*}
This is  a discretization of the following system of differential equations
\begin{eqnarray}
\begin{cases}
u'=2\big(1-v^{2}\big),\\
v'= 2\big(1-u^{2}\big).
\end{cases}
\label{1.4}\end{eqnarray}

We notice easily that the set $\mathfrak{E}=\big\{(u,v)\backslash
(u-v)(u^{2}+uv+v^{2}-3)=0 \big\}$ is an invariant set for this
system. The system (\ref{1.4}) has two  saddle points in $A=(1,1)$ and
$B=(-1,-1)$ and corresponding stable and unstable manifolds lie on its 
invariant set $\mathfrak{E}$.

The stable manifold at $A$ and unstable
manifolds at $B$  are part of set $\big\{(u,v)\backslash
u=v\big\}$, where as the unstable manifold at $A$ and stable
manifolds at $B$ are part of set  $\big\{(u,v)\backslash
u^{2}+uv+v^{2}=3\big\}$. The stable manifold of $A$ coincides with
 the unstable manifold of $B$ \,(See FIG.2).

 For the discretized equation (\ref{1.3}), these manifolds still exist \cite{FS}, let   
$W_{s}^{+}$ ,  $W_{i}^{-}$ denote the stable manifold  at $A$ and the
unstable manifold at $B$ respectively, $W_{i}^{+}$ and  $W_{s}^{-}$
unstable manifold at  $A$ and  stable manifold  at  $B$ respectively,
$W_{s}^{+}$ and  $W_{i}^{-}$ do not coincide an more (See FIG.3) as we 
want to show.

In the paper \cite{FS}, after introducing the notion of length of the 
first level $l_{1}(\varepsilon)=2\varepsilon n_{1}(\varepsilon)$, where 
$ n_{1}(\varepsilon)=\inf\big\{n\in\mathrm{IN}\backslash
y_{2n}(\varepsilon)+y_{2n+1}(\varepsilon)<0\big\}$, Fruchard-Sch\"afke
had shown that there exist a postive constant K such that 
\begin{eqnarray*}
\Delta(\varepsilon)&\leq&
\exp\bigg(-\frac{\pi^{2}+o(1)}{2\varepsilon}\bigg),\quad \text{as} \,\,\,
\varepsilon\searrow 0,\\
l_{1}(\varepsilon)&\geq& \frac{\pi^{2}}{4\varepsilon}+o(1),
\end{eqnarray*}
where $\Delta(\varepsilon)$ denote the distance between the sets
$W^{+}_{s}\cap S$ and $W^{-}_{i}\cap S$ in the sense of Hausdorff and
$S=\big\{(u,v)\backslash -1\leq u+v\leq 1\big\} $.

 With the initial condition $y_{0}=0,\,\, y_{1}=\varepsilon$, the length of the first level satisfies $l_{1}(\varepsilon)= 
-\frac{1}{2}\big(1+o(1)\big)\log(\Delta(\varepsilon))$ as
$\varepsilon\searrow 0$\cite{FS}. They also showed that there are two families of entire functions $y^{\pm}_{\varepsilon}:\mathbb{C}\longmapsto\mathbb{C}, $
solutions of the difference equation 
\begin{equation}
y(t+\varepsilon)=y(t-\varepsilon)+2\varepsilon\big(1-y(t)^{2}\big),
\label{1.5}\end{equation}
and the functions $t\mapsto
\big(y^{+}_{\varepsilon}(t),y^{+}_{\varepsilon}(t+\varepsilon)\big) $, 
$t\mapsto\big(y^{-}_{\varepsilon}(t),y^{-}_{\varepsilon}(t+\varepsilon)\big)
$ provide parametrization $w^{+}_{s}(t)$  of $W^{+}_{s}$ for $t\in [-1,\infty[$
respectively  $w^{-}_{i}(t)$  of $W^{-}_{i}$ for $t\in ]-\infty,1]$.

In this work we will prove 
\begin{theo}
There exist a constant $\alpha $ with $1.2641497\leq \alpha \leq
1.2641509$ and $\varepsilon_{0}>0$
such that  for $ 0<\varepsilon<\varepsilon_{0}$ 
\begin{equation*} 
Dist_{\varepsilon}\big(w^{+}_{s}(t),W^{-}_{i}\big)=\frac{4\pi\alpha\cos(\frac{\pi}{\varepsilon}t+\pi)}{\varepsilon^{3}\big(1-\tanh(t)^{2}\big)}e^{-\frac{\pi^{2}}{2\varepsilon}}+\mathcal{O}\Bigg(\frac{1}{\varepsilon^{2}}e^{-\frac{\pi^{2}}{2\varepsilon}}\Bigg),\quad \text{as} \,\,\,\varepsilon\searrow 0,
\end{equation*}
where 
$Dist_{\varepsilon}$ is defined as the vertical distance between the stable and
unstable manifolds. 
\label{t1.1}\end{theo}

 In order to show the theorem \ref{t1.1}, we start  with the construction
 of a formal power series solution   in $d$
 whose coefficients are polynomials in $u=\tanh(d\,t/\varepsilon),$  afterwards  we will
 give  asymptotic approximations  of these coefficients using appropriate 
 norms on the spaces of polynomials. The next step is to construct
 a quasi-solution i.e.\ a function that satisfies the equation (\ref{1.5})
 except for  an
 exponentially small error, then we show that this quasi-solution and the
 exact solution  of  equation (\ref{1.5})  are exponentially close.

 Finally we give an  asymptotic estimate of the  distance between the stable manifold of  $A$
 and the unstable manifold of $B$  and we show that this distance is
 exponentially small but not zero, thus completing the proof of the theorem \ref{t1.1}.

\begin{figure}[!t]
  \centering
  
   
  \includegraphics[width=.4\linewidth]{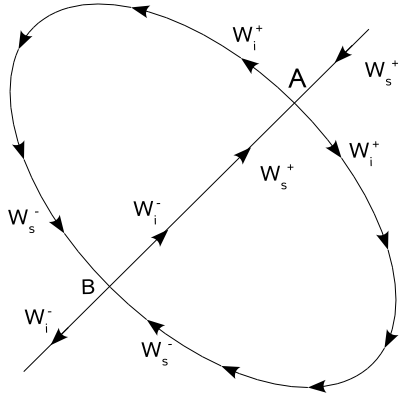}
  
  \caption{\small{The stable and unstable
      manifolds for the logistic differential  equation.}}
\end{figure}

\begin{figure}[!h]
  \centering
  
   
  \includegraphics[width=.4\linewidth]{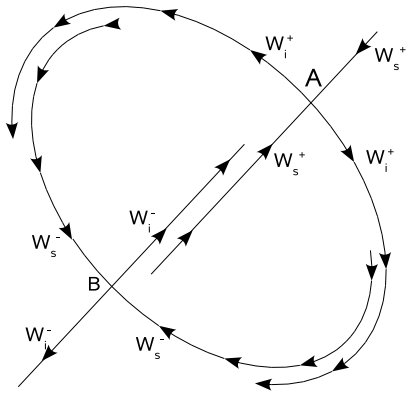}
  
  \caption{\small{The stable and unstable
      manifolds for the difference equation.}}
\end{figure}

\section{Formal solutions}

We consider the difference equation
\begin{equation}
y(t+\varepsilon)=y(t-\varepsilon)+2\varepsilon\big(1-y(t)^{2}\big),
\label{2.1}\end{equation}
where $\varepsilon>0 $ is the discretization step and $
  y(t) \to 1$ as $t \to +\infty$.
  Our first aim is to transform (\ref{2.1}) in such a way that the new equation
admits a formal solutions whose coefficients are polynomials.
 We define   $ u $ and $A(d,u)$ such that
 \begin{equation}
 u:=\tanh(\frac{d}{\varepsilon}t ),
\label{2.2}\end{equation}
\begin{equation}
y_{\varepsilon}(t)=A(d,u),
\label{2.3}\end{equation}
where $d:=\sum^{+\infty}_{n=1}{d_{n}\varepsilon^{n}}$ is a series to be determined.
When these variables are substituted in $(\ref{2.1})$, the following  equation satisfied
 by $ A(d,u)$ is obtained:
\begin{equation} 
A\big(d,T^{+}(d,u)\big)-A\big(d,T^{-}(d,u)\big)=f\big(\varepsilon,A(d,u)\big),
\label{2.4}\end{equation}
where
\begin{equation}
f(\varepsilon,x):=2\varepsilon(1-x^{2}),
\label{2.5}\end{equation}
\begin{equation}
T^{+}:=T^{+}(d,u)=\frac{u+\tanh(d)}{1+u\tanh(d)},
\label{2.6}\end{equation}
\begin{equation}
T^{-}:=T^{-}(d,u)=\frac{u-\tanh(d)}{1-u\tanh(d)}.
 \label{2.7}\end{equation}

 For small $\varepsilon $ one can construct a formal
expansion in powers of $\varepsilon^{2}$ of the form \quad
$\sum^{\infty}_{n=0} A_{2n+1}(u)d^{2n}$, where $A_{2n+1}$ are polynomials 
which satisfy  $A_{1}(u)=u$ and $A_{2n+1}(\pm 1)=0,$ for $n\geq 1$.

The existence of a such formal solution  is only possible if $d$ and 
 $\varepsilon$ are coupled in a very special way.  Indeed, 
  suppose that there exists  a such formal power series solution.
We differentiate (\ref{2.4}) with respect to u and obtain at $u=1$
\begin{displaymath} 
\frac{\partial A}{\partial u}(d,1)\cdot\lim_{u\to 1}\frac{T^{+}(d,u)-1}{u-1}-\frac{\partial A}{\partial u}(d,1)\cdot\lim_{u\to 1}\frac{T^{-}(d,u)-1}{u-1}=\frac{\partial f}{\partial x}(\varepsilon,1)\cdot\frac{\partial A}{\partial u}(d,1).
\end{displaymath}
 Because  $A_{2n+1}(1)=0$ for $n\geq 1$ implies $A(d,1)=1$. Thus with
$\frac{\partial f}{\partial x}(\varepsilon,1)=-4\varepsilon$ and  $\frac{\partial A}{\partial u}(d,1)=1+O(d)\neq 0$,  the above equation implies
 \begin{equation}
\varepsilon=\frac{\tanh(d)}{1-\tanh(d)^{2}}=\frac{1}{2}\sinh(2d).
 \label{2.8}\end{equation}

\begin{theo}
Let $\varepsilon$ and $d$ be coupled
by  $\varepsilon=\frac{1}{2}\sinh(2d)$. Then (\ref{2.4}) has a formal  power
series solution that can be written in the form \\

 \begin{equation}
A(d,u)= u+\sum^{\infty}_{n=1} A_{2n+1}(u)d^{2n},
\label{2.9}\end{equation}
where $A_{2n+1}(u)$ are odd polynomial and $A_{2n+1}(1)=A_{2n+1}(-1)=0$ for all $n$.
\label{t2.1}\end{theo} 

\noindent\textbf{Proof.}
 We will use the Induction Principle for even n  to show that there
 exist unique odd polynomials $A_{1},A_{3},A_{5}...A_{n+1}$ such that 
  \begin{equation}
Z(d,u)=\sum_{\substack{k=0\\ k\  even}}^{n}A_{k+1}d^{k}, 
\label{2.10}\end{equation} 
satisfy
\begin{equation}
Z\big(d,T^{+}(d,u)\big)-Z\big(d,T^{-}(d,u)\big)=f\big(\varepsilon,Z(d,u)\big)\
mod\  d^{n+2}.
\label{2.11}\end{equation}
\indent For $n=0$, we put $A_{1}(u)=u $ and $Z(d,u)=u$ and obtain 
\begin{displaymath}
Z\big(d,T^{+}(d,u)\big)-Z\big(d,T^{-}(d,u)\big)=T^{+}(d,u)-T^{-}(d,u)=(2-2u^2)d+\mathcal{O}(d^3).
\end{displaymath}
and

\begin{displaymath}
f\big(\varepsilon,Z(d,u)\big)=2\varepsilon(1-u^{2})=(2-2u^2)d+\mathcal{O}(d^3).
\end{displaymath}
This gives
\begin{displaymath}
Z\big(d,T^{+}(d,u)\big)-Z\big(d,T^{-}(d,u)\big)=f\big(\varepsilon,Z(d,u)\big)\ 
mod \  d^{2}.
\end{displaymath}
\indent Now suppose that for some even $n$ already
$A_{1},A_{3},A_{5}...A_{n+1}$ have been found with the above
properties. We have to construct $A_{n+3}$. First, we show that
$Z(d,u)$ satisfies (\ref{2.11}) even modulo $d^{n+3}$. To this purpose, let
\begin{equation}
Z\big(d,T^{+}\big)-Z\big(d,T^{-}\big)=f\big(\varepsilon,Z(d,u)\big)+ \mathit{R}_{n+2}(u)d^{n+2}+\mathcal{O}(d^{n+3})
\label{2.12}\end{equation} 
  we replace $d$ by $-d$. Using that $Z$ is   even and  $\varepsilon$ is odd  in $d$, we obtain 
\begin{displaymath}
Z\big(d,T^{-}\big)-Z\big(d,T^{+}\big)=-2\varepsilon\big(1-Z(d,u)^{2}\big)+ \mathit{R}_{n+2}(u)d^{n+2}+\mathcal{O}(d^{n+3}),
\end{displaymath}
this gives
\begin{displaymath}
Z\big(d,T^{+}\big)-Z\big(d,T^{-}\big)=2\varepsilon(d)\big(1-Z(d,u)^{2}\big)- \mathit{R}_{n+2}(u)(d)^{n+2}+\mathcal{O}(d^{n+3}).
\end{displaymath}
 With (\ref{2.12}) this  implies $\mathit{R}_{n+2}(u)=0$ and consequently
\begin{equation}
Z\big(d,T^{+}\big)-Z\big(d,T^{-}\big)=f\big(\varepsilon,Z(d,u)\big)+ \mathit{R}_{n+3}(u)d^{n+3}+\mathcal{O}(d^{n+4}).
\label{2.13}\end{equation} 
\indent We want to construct $A_{n+3}(u)$ such that
\begin{equation}
\widetilde{Z}\big(d,T^{+}\big)-\widetilde{Z}\big(d,T^{-}\big)=f\big(\varepsilon,\widetilde{Z}(d,u)\big)+
\mathcal{O}(d^{n+4}), 
\label{2.14}\end{equation}
if we put
\begin{equation} 
\widetilde{Z} (d,u) =Z(d,u)+A_{n+3}(u)d^{n+2}.
\label{2.15}\end{equation}
To this purpose  we use again Taylor expansion 
 \begin{eqnarray*}
 f\big(\varepsilon,\widetilde{Z} (d,u)\big) &=& f\big(\varepsilon,
 Z(d,u)\big)-4u\,A_{n+3}(u)d^{n+3}+\qquad\qquad\qquad\quad\  \mathcal{O}(d^{n+4}),   \\
 \widetilde{Z} \big(d,T^{+}(d,u)\big) &=& 
 Z(d,u)+A_{n+3}(u)d^{n+2}+(1-u^{2})\frac{\partial{A_{n+3}}}{\partial{u}}(u)d^{n+3}+\mathcal{O}(d^{n+4}),  \\
\widetilde{Z} \big(d,T^{-}(d,u)\big) &=& 
 Z(d,u)+A_{n+3}(u)d^{n+2}-(1-u^{2})\frac{\partial{A_{n+3}}}{\partial{u}}(u)d^{n+3}+\mathcal{O}(d^{n+4}).
  \end{eqnarray*}
  With (\ref{2.13}), this  gives
\begin{equation} 
\widetilde{Z}\big(d,T^{+}\big)-\widetilde{Z}\big(d,T^{-}\big)=f\big(\varepsilon,\widetilde{Z}(d,u)\big)+
\widetilde{R}_{n+3}(u)d^{n+3}+\mathcal{O}(d^{n+4}),
\label{2.16}\end{equation}
where
\begin{displaymath}
\widetilde{R}_{n+3}(u)=2(1-u^{2})\frac{\partial{A_{n+3}}}{\partial{u}}(u)+4u\,A_{n+3}(u)+\mathit{R}_{n+3}(u).
\end{displaymath}
We see that (\ref{2.14}) is satisfied if only if
\begin{equation} 
2(1-u^{2})\frac{\partial{A_{n+3}}}{\partial{u}}(u)+4u\,A_{n+3}(u)+\mathit{R}_{n+3}(u)=0.
\label{2.17}\end{equation}
This equation has a unique odd solution which is given by
\begin{eqnarray*}          
A_{n+3}(u)=-(1-u^{2})\int^{u}_{0}\frac{\mathit{R}_{n+3}(t)}{2(1-t^{2})^{2}}dt. 
\end{eqnarray*}

 We can prove that this solution is polynomial. For that, it is
 necessary that $\mathit{R}_{n+3}(u)$ and $\mathit{R}'_{n+3}(u)$ vanish into 1 and -1. 
Indeed, if we take  (\ref{2.13}) with $u=1$, we obtain 
\begin{equation}
\mathit{R}_{n+3}(1)d^{n+3}=Z\big(d,T^{+}(d,1)\big)-Z\big(d,T^{-}(d,1)\big)-f\big(\varepsilon,Z(d,1)\big)+\mathcal{O}(d^{n+4}).
\label{2.18}\end{equation}
Since,  $T^{+}(d,1)=T^{-}(d,1)=Z(d,1)=1$ and  $f\big(\varepsilon,1\big)=0$, 
we  obtain  
\begin{equation}
\mathit{R}_{n+3}(1)d^{n+3}=\mathcal{O}(d^{n+4}), 
\label{2.19}\end{equation}
  therefore $\mathit{R}_{n+3}(1)=0$.
 In order to show that $\mathit{R}'_{n+3}(1)=0$,  we derive formally (\ref{2.13}) and
 replace $u$ by 1. Using $T^{+}(d,1)=T^{-}(d,1)=Z(d,1)=1$ and  $\frac{\partial{f}}{\partial{u}}\big(\varepsilon,1\big)=-4\varepsilon$, we obtain 
\begin{equation}
\mathit{R}'_{n+3}(1)d^{n+3}=\frac{\partial{Z}}{\partial{u}}(d,1)\bigg(\frac{\partial{T^{+}}}{\partial{u}}(d,1)-\frac{\partial{T^{-}}}{\partial{u}}(d,1)+4\varepsilon\bigg)+\mathcal{O}(d^{n+4}).
\label{2.20}\end{equation}
By choice of $d$, cf (\ref{2.8}). Thus we have
$\frac{\partial{T^{+}}}{\partial{u}}(d,1)-\frac{\partial{T^{-}}}{\partial{u}}(d,1)=-2\sinh(2d)$. 
with (2,8) we obtain
\begin{equation}
\mathit{R}'_{n+3}(1)d^{n+3}=\mathcal{O}(d^{n+4}), 
\label{2.21}\end{equation}
therefore $\mathit{R}'_{n+3}(1)=0$. As formula (\ref{2.13}) shows that  $\mathit{R}_{n+3}(u)$ is odd, then 
$\mathit{R}_{n+3}(-1)=\mathit{R}'_{n+3}(-1)=0$.

The first polynomials  $A_{n+1}(u)$ with even $n$ are given by
$$\begin{array}{|c|c|c|c|c|} \hline

n&0 &2 &4     &   6\hl

 A_{n+1}(u)\  & u &u-u^{3} 
 &(1-u^{2})\Big(\frac{4}{3}u-\frac{10}{3}u^{3}\Big)&(1-u^{2})\Big(\frac{62}{3}u^{5}-\frac{190}{9}u^{3}+\frac{182}{45}u\Big)  \hl

\end{array}
$$

\indent We introduce the operators $\mathcal{C}_{2}, \mathcal{C},
\mathcal{S}_{2}, \mathcal{S}$ defined by

\begin{equation}          
\begin{array}{llll}
    \mathcal{C}(Z)(d,u)=\frac{1}{2}\big(Z(d,T^{+\frac{1}{2}})+Z(d,T^{-\frac{1}{2}})\big),&\\
\\
    \mathcal{S}(Z)(d,u)=\frac{1}{2}\big(Z(d,T^{+\frac{1}{2}})-Z(d,T^{-\frac{1}{2}})\big),&\\
\\
\mathcal{C}_{2}(Z)(d,u)=\frac{1}{2}\big(Z(d,T^{+})+Z(d,T^{-})\big),&\\
\\
\mathcal{S}_{2}(Z)(d,u)=\frac{1}{2}\big(Z(d,T^{+})-Z(d,T^{-})\big).&
\end{array}
\label{2.22}\end{equation}
where
$T^{+\frac{1}{2}}=T^{+}(\frac{d}{2},u)$,
$T^{-\frac{1}{2}}=T^{-}(\frac{d}{2},u)$ and $Z(d,u)$ is a  formal
power series of d whose coefficients are polynomials. We rewrite
equation (\ref{2.4}) as 
\begin{equation}
\mathcal{S}_{2}(A)(d,u)= \varepsilon \big( 1-A(d,u)^{2} \big).
\label{2.23}\end{equation}

\section{Norms for polynomials and basis }

In the sequel we denote: 
\begin{itemize}
\item $\mathcal{P}$ the set of all polynomial whose coefficents are complex,
\item $\mathcal{P}_{n}$ the spaces of all polynomials of degree less than or equal to n,
\item
  $\mathcal{Q}:=\left\{\mathnormal{Q}(d,u)=\sum^{\infty}_{n=0}\mathnormal{Q}_{n}(u)d^{n},\textrm{ where}\  \mathnormal{Q}_{n}(u)\in \mathcal{P}_{n},\textrm{ for all}\  n\in \mathbb{N} \right\}.$
\end{itemize}
\begin{propo} \cite{SV}
If  we define the sequence of the polynomial functions $ \tau_{n}(u)$
by $\tau_{0}(u)=1, \tau_{1}(u)=u, \tau_{n+1}(u)=\frac{1}{n}D\tau_{n}(u)$, where the operator $D$ is defined by \quad $D:=(1-u^{2})\frac{\partial}{\partial{u}}$, we have
\begin{enumerate}
\item $T^{+}(d,u)=\sum^{\infty}_{n=0}\tau_{n+1}(u)d^{n},$
\item $ \tau_{n}(u)$   has exactly
 degree $n$,
\item $\tau_{n}(\tanh(z))=\frac{1}{(n-1)\,!}\big(\frac{d}{dz} \big)^{n-1}\big(\tanh(z)\big).$
\end{enumerate}
\label{p3.1}\end{propo}
\begin{defi} Let $p \in \mathcal{P}_{n}$.  As $
  \tau_{0}(u), \tau_{1}(u),...,\tau_{n}(u)$ form a basis of $\mathcal{P}_{n}$,   we can write
  $p$ as 
\begin{eqnarray*}
 p=\sum^{n}_{k=0}a_{k}\tau_{k}(u).
\end{eqnarray*}
Then we define the norm
 \begin{equation}
\|p\|_{n}=\sum^{n}_{i=0}a_{i}\left(\frac{\pi}{2}\right)^{n-i} .
\label{3.1}\end{equation}
\label{d3.2}\end{defi}\

\begin{theo}\cite{SV} Let n,m be positive integers and $p \in \mathcal{P}_{n}$, 
  \,$q \in \mathcal{P}_{m}$.
  The norms $\|\ \|_{n}$ of the above definition   have the following property:
\begin{enumerate}

\item  $ \|Dp\|_{n+1}\le n\|p\|_{n}$.   
\item If $p$ odd we have $\|p\|_{n}\le \|Dp\|_{n+1}$.
\item There exists a constant $M_{1}$ such that $\|pq\|_{n+m}\le M_{1} \|p\|_{n}|q\|_{m}$.
\item There is a constant $M_{2}$ such that $|p(u)|\le
  M_{2}\left(\frac{2}{\pi}\right)^{n} \|p\|_{n}$ $(-1\le u\le 1)$.

\item  There is a constant $M_{3}$ such that for all $ n> 1$  with
  $p(-1)= p(-1)=0,$
\begin{eqnarray*}
\Big\|\frac{p}{\tau_{2}}\Big\|_{n-2}\le M_{3} \|p\|_{n}.
\end{eqnarray*}
\end{enumerate}
\label{t3.3}\end{theo}
\section{Operators}
In this section we will use definitions and results adapted from
\cite{SV} by replacing $\frac{\pi}{2}$ by  $\pi$.
\begin{defi}
Let $f$ be  formal power series of $z$ whose coefficients are complex. 
We define a linear operator $f(dD)$ on $\mathcal{Q}$ by
\begin{equation}
f(dD)\mathnormal{Q}=\sum^{\infty}_{n=0} \Big( \sum^{n}_{i=0} f_{i}
   D^{i}\mathnormal{Q}_{n-i}(u) \Big)d^{n},
\label{4.1}\end{equation}
Where $f(z)=\sum^{\infty}_{i=0}f_{i}z^{i}$ and $ \mathnormal{Q}=\sum^{\infty}_{n=0} \mathnormal{Q}_{n}(u)d^{n}\in\mathcal{Q}$.
\label{d4.1}\end{defi}
 By the above definition and (1) of  the proposition \ref{p3.1}, we can show that 

\begin{eqnarray} 
\mathnormal{Q}(d,T^{+}\big(\theta d, u)\big)=\big(\exp(\theta
dD)\mathnormal{Q}\big)(d,u), \quad \textrm{for }\ \mathnormal{Q}\in \mathcal{Q}
 \  \textrm{and all}\  \theta\in \mathbb{C}.
\label{4.2}\end{eqnarray} 
With (\ref{2.22})  this implies
\begin{equation*}          
\begin{array}{llll}
    \mathcal{C}(\mathnormal{Q})(d,u)\
    =\cosh(\frac{d}{2}D)\mathnormal{Q},\ \mathcal{S}(\mathnormal{Q})(d,u)\ =\sinh(\frac{d}{2}D)\mathnormal{Q},&\\ 
\mathcal{C}_{2}(\mathnormal{Q})(d,u)=\cosh(dD)\mathnormal{Q},\
\mathcal{S}_{2}(\mathnormal{Q})(d,u)=\sinh(dD)\mathnormal{Q},
\end{array}
\end{equation*}
for polynomial series  $\mathnormal{Q}$ in $ \mathcal{Q}$. We denote $\big\|\mathnormal{Q}\big\|_{n}=\big\|\mathnormal{Q}_{n}\big\|_{n}$.
\begin{theo}.
Let $f(z)$ be  formal power series having a radius of convergence
greater than $\pi$ and let k be a positive integer. There is a constant 
$K$ such that: 
If $\mathnormal{Q}$ is a polynomial series having the following property 

\begin{eqnarray*} 
\|\mathnormal{Q}\|_{n}\leq \left\{\begin{array}{ll} 0 & \textrm{for}\  n<k
      \\ M(n-k)\,!\pi^{-n}  & \textrm{for}\  n\geq k
 \end{array} \right.
\end{eqnarray*} 
where $M$ independent of $n$ and  $\mathnormal{Q} \in  \mathcal{Q}$
then the polynomial series $f(dD)\mathnormal{Q}$ satisfies 
\begin{eqnarray*} 
\|f(dD)\mathnormal{Q}\|_{n}\leq \left\{\begin{array}{ll} 0 & \textrm{for}\  n<k
      \\ MK(n-k)\,!\pi^{-n}  & \textrm{for}\  n\geq k
 \end{array} \right.
\end{eqnarray*} 
\label{t4.2}\end{theo}
We define the operator $\mathcal{C}^{-1}$ by $\mathcal{C}^{-1}=g(dD)$, where $g(z)=\dfrac{1}{\cosh(z/2)}$
\begin{theo}   
There exists a positive constant $K$
 such that, if $\mathnormal{Q}$ is a polynomial series 
such that $\mathnormal{Q}_{n}$ odd, $\|\mathnormal{Q}\|_{n}=0$ for
$n<k$ for some positive integer $k$ and  
\begin{eqnarray*} 
\|dD\mathnormal{Q}\|_{n}\leq M(n-k)\,!\pi^{-n}  & \textrm{for}\  n\geq k
\end{eqnarray*} 
 where $M$ independent of $n$ and  $\mathnormal{Q} \in  \mathcal{Q}$, the polynomial series
 $\mathcal{C}^{-1}(\mathnormal{Q})$ satisfies 
\begin{eqnarray*} 
\big\|\mathcal{C}^{-1}(\mathnormal{Q})\big\|_{n}\leq  MK\pi^{-n} \left\{\begin{array}{ll}
 n\,! & \textrm{for}\  k=1
 \\ (n-1)\,!\log(n)  & \textrm{for}\   k=2 
\\ (n-1)\,! & \textrm{for}\ k \geq 3 
 \end{array} \right.
\end{eqnarray*}
\label{t4.3}\end{theo}
\begin{theo}
We consider a polynomial series 
\begin{eqnarray*} 
\mathnormal{Q}_{\alpha}(d,u)=\sum^{\infty}_{\substack{n=1\\ n\  odd} }\alpha_{n}(n-1)\,!\Big(
\frac{i}{\pi}\Big)^{n-1}\tau_{n}(u)d^{n}
\end{eqnarray*}
where $\alpha_{n}=\mathcal{O}(n^{-k})$ as $n
\to \infty$ with some integer $k\geq 2$. Let $\alpha:=
\frac{1}{\pi}\sum^{\infty}_{n=1} \alpha_{n}$. then the coefficients  
$\left\{ \mathcal{C}^{-1}(\mathnormal{Q}_{\alpha}) \right\}_{n}$ of\,
$\mathcal{C}^{-1}(\mathnormal{Q}_{\alpha}) $ satisfy
\begin{eqnarray*} 
\Bigg\| \left\{ \mathcal{C}^{-1}(\mathnormal{Q}_{\alpha})
\right\}_{n}-\alpha (n-1)\,!\Big(\frac{i}{\pi}\Big)^{n-1}\tau_{n}\Bigg\|_{n}=\mathcal{O}\Big((n-k)\,!\pi^{-n}\Big),
\end{eqnarray*}
as $n\to \infty$ for odd $n.$
\label{t4.4}\end{theo}
\begin{theo}
Let $k,l, p,q$ be positive integer with $p\geq k$ and $q\geq l$. Define
$m$ as the minimum of $k+q$ and $l+p$. Then there is a constant $K$
with the following property:\\
\indent If $\mathnormal{P}$ and $\mathnormal{Q}$ are polynomial series 
such that $\|\mathnormal{P}\|_{n}=0 $ for
$n<p$,\ $\|\mathnormal{Q}\|_{n}=0 $ for $n<q$ and 
\begin{eqnarray*}
\begin{array}{ll}
 \|\mathnormal{P}\|_{n}\leq M_{1}(n-k)\,!\, \pi^{-n},  &
 \textrm{for}\,\,  n\geq k, \\ 
\|\mathnormal{Q}\|_{n}\leq  M_{2}(n-l)\,!\, \pi^{-n},  & \textrm{for}\,\,   
 n\geq l.   
\end{array}
\end{eqnarray*}
Then
\begin{eqnarray*}
\|\mathnormal{PQ}\|_{n}\leq KM_{1} M_{2}(n-m)\,!\, \pi^{-n},  & \textrm{for}\,\,   
 n\geq p+q.
\end{eqnarray*}
\label{t4.5}\end{theo}
\begin{theo}
Let $\mathnormal{Q}_{1}(d,u)$ be a convergent polynomial series which
is even with respect to both variables and has constant term 1.\\
 Let $\mathnormal{Q}_{2}(d,u)=d^{2}(1-u^{2})\mathnormal{Q}_{1}(d,u)$
 and
 $\mathnormal{P}(d,u)=\mathcal{S}(\mathnormal{Q}_{2})/\mathcal{C}(\mathnormal{Q}_{2})$. Consider the linear operator defined by 
\begin{equation}
\mathcal{L}(\mathnormal{Q})=\mathcal{S}(\mathnormal{Q})-\mathnormal{P}(d,u)\cdot\mathcal{C}(\mathnormal{Q}),\quad 
\mathnormal{Q}\in \mathcal{Q}.
\label{4.3}\end{equation}
Then, there exist a constant $K$ with the following property. If
$\mathnormal{Q}$ is an odd polynomial series with odd coefficients
$\mathnormal{Q}_{n}(u)$ satisfying $\mathnormal{Q}_{n}(1)=0$ for
all\ $n$, $\|\mathcal{L}(\mathnormal{Q})\|_{n}=0$ for $n<6$ and 
\begin{eqnarray*}
\|\mathcal{L}(\mathnormal{Q})\|_{n}\leq M (n-6)\,!\pi^{-n}\  \textrm{for}\ 
n\geq 6,
\end{eqnarray*}
then also 
\begin{eqnarray*}
\|dD\mathnormal{Q}\|_{n}\leq K M (n-6)\,!\pi^{-n}\  \textrm{for}\ 
n\geq 6.
\end{eqnarray*}
\label{t4.6}\end{theo}

\section{Asymptotic approximation of the coefficients of the formal solution} 

The objective in this section is to construct an asymptotic aproximation
 of the coefficents of the formal solution (\ref{2.9}).  It will
turn out to be convenient to consider the new series $ B(d,u)=\varepsilon A(d,u), $ 
 \begin{eqnarray*}
B(d,u)=\varepsilon A(d,u)=\sum^{\infty}_{\substack{n=2\\ n\
    even}}A_{n+1}(u)\varepsilon d^{n}+\varepsilon u,
\end{eqnarray*}
this gives
\begin{eqnarray*}
B(d,u)=\sum^{\infty}_{\substack{n=3\\ n\
    odd}}B_{n}(u) d^{n}+\varepsilon u,
\end{eqnarray*}
where $B_{n}(u) $ are odd polynomials. Furthermore we have
$B_{n}(u)\in \mathcal{P}_{n}$ and $B_{n}(1)= B_{n}(-1) =0 $ for all $n$. The new equation for  $B(d,u)$
is
\begin{eqnarray}
\mathcal{S}_{2}(B)(d,u)=\varepsilon^{2} -B(d,u)^{2}.
\label{5.1}\end{eqnarray}
We saw in section 2 that the series $A$ was a formal solution of the
starting equation. For the moment, nothing is known about the norms of its coefficients, but we will show that this series is Gevrey-1, more precisely $\|A\|_{n}=O\big(n\,!\,\pi^{-n}\big)$. 
 This enables us thereafter to construct a quasi-solution.  To this purpose, we will prepare the equation (\ref{5.1}) so that we can construct a recurrence.
 
 We start with the decompostion of series $B(d, u)$ in the form:
\begin{eqnarray}
B=U+F, \,\textrm{ where}\ \ \  U=\varepsilon u +(u-u^{3})d^{3}+\Big(\frac{10}{3}u^{5}-\frac{16}{3}u^{3}+2u\Big)d^{5}.
\label{5.2}\end{eqnarray}
 Next, we define
\begin{eqnarray}
G:&=& Q\,F,\nonumber \\
J:&=&\frac{\mathcal{C}(G)}{Q_{1}}.
\label{5.3}\end{eqnarray}
where 
\begin{eqnarray}
 Q(d,u):&=& (1-u^{2})d^{2}
  +(u^{4}-u^{2})d^{4}-\Big(\frac{13}{6}u^{6}-\frac{7}{2}u^{4}+\frac{4}{3}u^{2}\Big)d^{6}\nonumber\\
&+&\Big(\frac{47}{6}u^{8}-\frac{31}{6}u^{6}+\frac{58}{15}u^{4}-\frac{104}{45}u^{2}+1\Big)d^{8},\nonumber\\
Q_{1}(d,u)&=&(1-u^{2})d^{2}+\frac{3}{2}(u^{2}-u^{4})d^{4}.
\label{5.4}\end{eqnarray}

The choice of $Q$ and $Q_{1}$ and by using the properties of the operators  $\mathcal{S},\,  \mathcal{C}$,  we will be able to rewrite the equation (\ref{5.1}) in the form:
\begin{eqnarray*}
e_{0}(d,u)\mathcal{L}(J)&=&e_{1}(d,u)\mathcal{C}(J)+ e_{2}(d,u)F+e_{3}(d,u)F^{2} \nonumber\\ & &
 +e_{4}(d,u)\mathcal{C}_{2}(F)+e_{5}(d,u)  
\end{eqnarray*}
  where  $\mathcal{L}$ is the operator defined in (\ref{4.3}) and
  $e_{i}(d,u),i=0,...,5 $  are convenient known convergent series in
  $d$ and $u$ which will be thereafter given.

 The left hand side of this equation is series 
  with leading term 1 multiplied by the invertible operator $\mathcal{L}$ applied to the
 series $J$. The right hand side  is an expression of $F$  and $J$
 multiplied by known convergent  series  $e_{i}(d,u), i=1,...,5$.

 $U$, $Q$, $Q_{1}$   were chosen so that  the series
$e_{i}(d,u), i=0..5$ are of a rather large order in $d$, this makes
the second term smaller than the right hand side. This property will be useful to construct a recurrence on $n$ and to
 reverse then the operators $\mathcal{L}, \mathcal{C}$, which makes
 possible to estimate the coefficients of series $\mathcal{C}(G)$ and
 finally  the coefficients of the formal solution of the
 equation (\ref{5.1}). 

More precisely, we insert (\ref{5.2}) into the equation (\ref{5.1}) and find 
\begin{equation}
\mathcal{S}_{2}(F)= -F^{2}-2U\,F +\varepsilon^{2}-U^{2}-\mathcal{S}_{2}(U).
\label{5.5}\end{equation}
We define  $\mathcal{X}$ by
\begin{eqnarray}
\mathcal{X}:= \mathnormal{V}(d,u)\cdot\mathcal{S}\mathcal{C}(G)+
\mathnormal{W}(d,u)\cdot\mathcal{C}^{2}(G)
\label{5.6}\end{eqnarray}
where 
\begin{eqnarray*}
 \mathnormal{V}(d,u):&=&1+2u^{2}d^{2}-\big(\frac{7}{3}u^{4}-3u^{2}\big)d^{4}-\Big(\frac{31}{3}u^{6}-\frac{221}{9}u^{4}+\frac{253}{15}u^{2}-\frac{41}{15}\Big)d^{6}, \nonumber\\
\mathnormal{W}(d,u):&=&2\varepsilon u. 
\end{eqnarray*}
Then, because of   $G=Q\,F $, 
we have \  $\mathcal{X}(G)= \mathnormal{V}\cdot\mathcal{S}\mathcal{C}(F\,Q)+
\mathnormal{W}\cdot\mathcal{C}^{2}(F\,Q).$
  Using the  formulas 
\begin{equation}
\begin{array}{llll}
   \mathcal{C}_{2}=2\mathcal{C}^{2}-Id,  &\\
\\
 \mathcal{S}_{2}=2\mathcal{S}\mathcal{C},   &\\
\\
 \mathcal{S}_{2}(F\,Q)= \mathcal{S}_{2}(F)\mathcal{C}_{2}(Q)+\mathcal{S}_{2}(Q)\mathcal{C}_{2}(F),  &\\
\\
\mathcal{C}_{2}(F\,Q)= \mathcal{C}_{2}(F)\mathcal{C}_{2}(Q)+\mathcal{S}_{2}(Q)\mathcal{S}_{2}(F). &
\end{array}
\label{5.7}\end{equation}
We obtain 
\begin{eqnarray*}
\mathcal{X}&=&\frac{1}{2}\mathnormal{V}\cdot\mathcal{S}_{2}(F\,Q)+\frac{1}{2}\mathnormal{W}\cdot\Big(\mathcal{C}_{2}(F\,Q)+Q\,F\Big) \\
      &=& \frac{1}{2}\Big( \mathnormal{V}\,\mathcal{C}_{2}(Q)+\mathnormal{W}\,\mathcal{S}_{2}(
Q)\Big)\mathcal{S}_{2}(F)
+\frac{1}{2}\Big(\mathnormal{V}\,\mathcal{S}_{2}(Q)+\mathnormal{W}\,\mathcal{C}_{2}(Q)\Big)\mathcal{C}_{2}(F)\nonumber\\ 
& & {}+\frac{1}{2}\mathnormal{W}\,Q\cdot F.  
\end{eqnarray*}
With (\ref{5.5}) this implies
\begin{eqnarray}
\mathcal{X}=  \frac{1}{2}\mathnormal{W}_{1}\,F
+\frac{1}{2}\mathnormal{W}_{2}\,\mathcal{C}_{2}(F)
-\frac{1}{2}\big(\mathnormal{V}\,\mathcal{C}_{2}(Q)+\mathnormal{W}\,\mathcal{S}_{2}(Q)\big) 
 F^{2} +
\mu(\varepsilon,u), 
\label{5.8}\end{eqnarray}
where 
\begin{eqnarray*}
\begin{array}{llll}
 \mathnormal{W}_{1}=-2U\big(\mathnormal{V}\,\mathcal{C}_{2}(Q)+\mathnormal{W}\,\mathcal{S}_{2}(Q)\big)+\mathnormal{W}\,Q,&\\
\\
\mathnormal{W}_{2}=\mathnormal{V}\,\mathcal{S}_{2}(Q)+\mathnormal{W}\,\mathcal{C}_{2}(Q),&\\ 
\\
\mu(\varepsilon,u)=\frac{1}{2}\big(\mathnormal{V}\,\mathcal{C}_{2}(Q)+\mathnormal{W}\,\mathcal{S}_{2}(Q)\big)\big(\varepsilon^{2}-U^{2}-\mathcal{S}_{2}(U)\big). &
\end{array}
\end{eqnarray*}

The lowest power of $d$ in $\mathnormal{W}_{1}$ and
$\mathnormal{W}_{2}$ is the eleventh 
$\mathnormal{W}_{1}=\mathcal{O}(d^{11})=\mathnormal{W}_{2}$, $
\mu(\varepsilon,u)$ is analytic. On the other hand $\mathcal{C}(G)=Q_{1}\,J,$
this implies
\begin{equation}
\mathcal{X}(G)= \mathnormal{V}\cdot\mathcal{S}(Q_{1}\,J)+
\mathnormal{W}\cdot\mathcal{C}(Q_{1}\,J).
\label{5.9}\end{equation}
If we use the product formulas
\begin{equation}
\begin{array}{ll}
  
 \mathcal{S}(Q_{1}\,J)= \mathcal{S}(J)\mathcal{C}(Q_{1})+\mathcal{S}(Q_{1})\mathcal{C}(J),  &\\
\\
 \mathcal{C}(Q_{1}\,J)= \mathcal{C}(J)\mathcal{C}(Q_{1})+\mathcal{S}(Q_{1})\mathcal{S}(J). &
\end{array}
\label{5.10}\end{equation}
We have 
\begin{equation}
\mathcal{X}(G)=\Big(\mathnormal{V}\mathcal{C}(Q_{1})+\mathnormal{W}\mathcal{S}(Q_{1})\Big)\mathcal{S}(J)+
\Big(\mathnormal{V}\mathcal{S}(Q_{1})+\mathnormal{W}\mathcal{C}(Q_{1})\Big)\mathcal{C}(J).
\label{5.11}\end{equation}
With (\ref{5.8}) and (\ref{4.3}), this implies 
\begin{eqnarray}
\mathnormal{V}_{1}\mathcal{L}(J)&=&-(\mathnormal{W}_{3}+P\mathnormal{V}_{1})\mathcal{C}(J)+ \frac{1}{2}\mathnormal{W}_{1}F\,
-\frac{1}{2}\big(\mathnormal{V}\,\mathcal{C}_{2}(Q)+\mathnormal{W}\,\mathcal{S}_{2}(Q)\big) 
 F^{2} \nonumber\\ & &
 {}+\frac{1}{2}\mathnormal{W}_{2}\,\mathcal{C}_{2}(F)
 +\mu(\varepsilon,u), 
\label{5.12}\end{eqnarray} 
\begin{eqnarray}
 \mathnormal{V}_{1}&=& \mathnormal{V}\mathcal{C}(Q_{1})+\mathnormal{W}\mathcal{S}(Q_{1}), \nonumber \\
 \mathnormal{W}_{3}&=&\mathnormal{V}\mathcal{S}(Q_{1})+\mathnormal{W}\mathcal{C}(Q_{1}), \nonumber\\
    P&=&
    \mathcal{S}\big((1-u^{2})d^{2}\big)/\mathcal{C}\big((1-u^{2})d^{2}\big)=-ud +\mathcal{O}(d^{2}), \nonumber \\
\mathcal{L}(J)&=&\mathcal{S}(J)-P\cdot\mathcal{C}(J).
\label{5.13}\end{eqnarray}

We divide (\ref{5.12})  by $(1-u^{2})d^{2}$, this implies
\begin{eqnarray}
\frac{\mathnormal{V}_{1}}{(1-u^{2})d^{2}}\mathcal{L}(J)&=&-\frac{\mathnormal{W}_{3}+P\mathnormal{V}_{1}}{(1-u^{2})d^{2}}\mathcal{C}(J)+ \frac{\mathnormal{W}_{1}}{2(1-u^{2})d^{2}}F
-\frac{\mathnormal{V}\,\mathcal{C}_{2}(Q)+\mathnormal{W}\,\mathcal{S}_{2}(Q)}{2(1-u^{2})d^{2}} 
 F^{2} \nonumber\\ & &
 {}+\frac{\mathnormal{W}_{2}}{2(1-u^{2})d^{2}}\mathcal{C}_{2}(F)+\frac{\mu(\varepsilon,u)}{(1-u^{2})d^{2}}, 
\label{5.14}\end{eqnarray}
where $F$ and $J$  are coupled by the equation $Q_{1}J=\mathcal{C}(Q\,F)=\mathcal{C}(G).$
\begin{theo}
With the above notations, we have the following estimation
\begin{displaymath}
\big\|dD\mathcal{C}(G)\big\|_{n+1}=\mathcal{O}\Big( (n-7)\,!\pi^{-n}\Big)
\quad \textrm{as}\  n\to \infty.
\end{displaymath}
\label{t5.1}\end{theo}
\noindent\textbf{Proof.} We set
\begin{equation}
e_{n}:=\frac{\big\|dD\mathcal{C}(G)\big\|_{n+1}}{ (n-6)\,!\pi^{-n}}\ \ 
\textrm{for}\ \  n \geq 7. 
\label{5.15}\end{equation}

We must show that $e_{n}=\mathcal{O}(n^{-1})$. In the
sequel we will use the convention: If  $a_{n}, n=0, 1, ...$. is any
sequence of positive real numbers, then 
\begin{equation}
a^{+}_{n}:=\max(a_{0},a_{1},...a_{n}), \ \textrm{for all} \ n\geq 0.
\label{5.16}\end{equation}

In this proof $K_{1}, K_{2},...$ will always denote constants
independent of $n$. Theorem \ref{t4.3} gives 
\begin{equation}
\big\|G\big\|_{n} \leq K_{1}e^{+}_{n}(n-1)\,!\pi^{-n},
\label{5.17}\end{equation}
\begin{equation}
\big\|G^{2}\big\|_{n} \leq K_{2}f_{n}(n-10)\,!\pi^{-n},\ \ \ \textrm{for}\ n\geq 18,
\label{5.18}\end{equation}
where 
\begin{displaymath}
f_{n}=
 \sum^{n-9}_{i=9}e^{+}_{i}e^{+}_{n-i} \frac{(i-1)\,!\,
  (n-i-1)\,!}{(n-10)\,!},\quad \textrm{for\ $n\geq 18$}.
\end{displaymath}
Using theorem \ref{t3.3}, (5), we find 
\begin{eqnarray}
\begin{array}{ll}
\big\|F\big\|_{n} \leq K_{4}e^{+}_{n+1}(n+1)\,!\pi^{-n},\\
\\
\big\|F^{2}\big\|_{n} \leq K_{3}f_{n+4}(n-6)\,!\pi^{-n}.
\end{array}
\label{5.19}\end{eqnarray}
Using theorem \ref{t4.2}, we obtain
\begin{eqnarray}
\big\|\mathcal{C}_{2}(F)\big\|_{n} \leq K_{5}e^{+}_{n+1}(n+1)\,!\pi^{-n}. 
\label{5.20}\end{eqnarray}

Observe that the convergent polynomials series 
$\mathnormal{W}_{1}/2(1-u^{2})d^{2}$\ ,$\ \mathnormal{W}_{2}/2(1-u^{2})d^{2}$ 
begins with $d^{9}$ and
$\big(\mathnormal{V}\,\mathcal{C}_{2}(Q)+\mathnormal{W}\,\mathcal{S}_{2}(Q)\big)/2(1-u^{2})d^{2}
$ begin with $1$.  Using   theorem \ref{t4.5} we thus obtain 
\begin{equation}
\Big\|\frac{\mathnormal{W}_{2}}{2(1-u^{2})d^{2}}\mathcal{C}_{2}(F)\Big\|_{n} \leq K_{6}e^{+}_{n-7}(n-7)\,!\pi^{-n},
 \label{5.21}\end{equation}
\begin{equation}
\Big\|\frac{\mathnormal{W}_{1}}{2(1-u^{2})d^{2}}F\Big\|_{n} \leq K_{7}e^{+}_{n-7}(n-7)\,!\pi^{-n},
 \label{5.22} \end{equation}
\begin{equation}
\Big\|\frac{\big(\mathnormal{V}\,\mathcal{C}_{2}(Q)+\mathnormal{W}\,\mathcal{S}_{2}(Q)\big)}{2(1-u^{2})d^{2}} F^{2}\Big\|_{n} \leq K_{8}f^{+}_{n+4}(n-6)\,!\pi^{-n}.
\label{5.23} \end{equation}
On the  other hand 
\begin{eqnarray*}
\big\|Q_{1}J\big\|_{n}=\big\| \mathcal{C}(G)\big\|_{n} \leq
e_{n}(n-6)\,!\pi^{-n}. 
\end{eqnarray*}
We apply the theorems \ref{t3.3}, (5) and \ref{t4.2} and find
\begin{eqnarray*}
\big\|\mathcal{C}(J)\big\|_{n} \leq K_{9}e_{n+2}(n-4)\,!\pi^{-n}. 
\end{eqnarray*}
Because the convergent polynomial series
$(\mathnormal{W}_{3}+P\mathnormal{V}_{1})/(1-u^{2})d^{2}=\mathcal{O}(d^{3}),$ 
then 
\begin{eqnarray}
\bigg\|\frac{(\mathnormal{W}_{3}+P\mathnormal{V}_{1})}{(1-u^{2})d^{2}}\mathcal{C}(J)\bigg\|_{n} \leq K_{10}e^{+}_{n-1}(n-9)\,!\pi^{-n}.  
\label{5.24}\end{eqnarray}
With a crude estimate of the convergent terms $\mu(\varepsilon,u)$ the 
inequalities (\ref{5.21})- (\ref{5.24}) gives
\begin{eqnarray*}
\Big\|\frac{\mathnormal{V}_{1}}{2(1-u^{2})d^{2}}\mathcal{L}(J)\Big\|_{n} \leq
K_{11}(1+e^{+}_{n-1}+f^{+}_{n+4})(n-6)\,!\pi^{-n}.
\end{eqnarray*}
The multiplication by the convergent terme
$2(1-u^{2})d^{2}/\,\mathnormal{V}_{1}$ only changes the constant
$K_{11}$. Since theorem 4.5 applies to the operator $\mathcal{L}$
defined in (\ref{4.3}), we obtain 
\begin{eqnarray}
\big\|dDJ\big\|_{n} \leq K_{12}(1+e^{+}_{n-1}+f^{+}_{n+4})(n-6)\,!\pi^{-n}.
\label{5.25}\end{eqnarray}
This with theorem \ref{t4.5} gives via 
\begin{eqnarray*}
dD\mathcal{C}(G)=dD(Q_{1} J)=J\cdot (dDQ_{1})+Q_{1}\cdot (dDJ),
\end{eqnarray*}
the estimate
\begin{eqnarray}
\big\|dD\mathcal{C}(G)\big\|_{n} \leq K_{13}(1+e^{+}_{n-3}+f^{+}_{n+2})(n-8)\,!\pi^{-n}.
\label{5.26}\end{eqnarray}
By equation (\ref{5.15}), we obtain
\begin{equation}
e_{n-1}\leq \frac{K}{n}(1+e^{+}_{n-3}+f^{+}_{n+2}).
\label{5.27}\end{equation}

\begin{lem} Under the condition (\ref{5.27}),  we have $e_{n}=\mathcal{O}(n^{-1})$ as
$n\to \infty .$
\end{lem} 
\noindent\textbf{Proof.} 
Let $K_{1}\geq 9\,!e^{+}_{10}$ an arbitrary number. We assume that 

\begin{equation}
e_{n}\leq \frac{K_{1}(n+p)\,!}{(n-1)\,!(p+10)\,!},\ \ \qquad  \textrm{for $10\leq  n\leq N-3 $},
\label{5.28}\end{equation}
with  $p\geq -1,\, N\geq 14$. This gives for $13\leq  n\leq N $
\begin{eqnarray*}
(n-8)\,!f_{n+2}\leq 2e^{+}_{9}K_{1}
\frac{ 8\,!(n+p-7)\,!}{!(p+10)\,!} +K_{1}^{2}\sum^{n-8}_{i=10}\frac{(i+p)\,!(n+p-i+2)\,!}{\big((p+10)\,!\big)^{2}}.
\end{eqnarray*} 
Using the inequality,
\begin{eqnarray*}
\sum^{n-8}_{i=10}(i+p)\,!(n+p-i+2)\,!\leq (p+10)\,!(n+p-7)\,!,
\end{eqnarray*} 
we obtain 
\begin{eqnarray*}
f_{n+2}\leq K_{2}\frac{(n+p-7)\,!}{(p+10)\,!(n-8)\,!}\leq K_{2}\frac{(n+p-3)\,!}{(p+10)\,!(n-4)\,!}  \qquad  \textrm{for $13\leq  n\leq N $}, 
\end{eqnarray*} 
 with a constant $K_{2}$ depends on $k_{1}$, independent of $p$.
The asumption (\ref{5.27}) of the lemma yields 
\begin{eqnarray*}
e_{n-1}\leq\frac{K}{n} (1+K_{1}+K_{2})\frac{(n+p-3)\,!}{(p+10)\,!(n-4)\,!}  \qquad  \textrm{for $13\leq  n\leq N $},  
\end{eqnarray*} 
this implies
\begin{eqnarray*}
e_{n}\leq\frac{K_{3}}{n+1}\frac{(n+p)\,!}{(p+10)\,!(n-1)\,!},  \qquad  \textrm{for $12\leq  n\leq N-1 $},  
\end{eqnarray*} 
with a constant $K_{3}$ only depending upon $K_{1}$.
Now we choose $N_{0}\geq 12$ so large that $\frac{K_{3}}{N_{0}}\leq
K_{1}$ and then $p$ so large that (\ref{5.28}) holds. Then 

\begin{eqnarray*}
e_{n}\leq\frac{K_{1}}{n+1}\frac{(n+p)\,!}{(p+10)\,!(n-1)\,!}  \qquad
\textrm{for $12\leq  n\leq N-1 $}, \, \textrm{with $ N\geq N_{0}$}.
\end{eqnarray*}
Since $K_{1}\geq 9\,! e_{10}$ is an arbitrary number, we have shown for 
any $p\geq -1$ , taht 
\begin{eqnarray*}
e_{n}=\mathcal{O}\bigg(\frac{(n+p)\,!}{(n-1)\,!}\bigg) \ \textrm{as $  n\to \infty $},
\end{eqnarray*}
implies that
\begin{eqnarray*}
e_{n}=\mathcal{O}\bigg(\frac{(n+p-1)\,!}{(n-1)\,!}\bigg) \ \textrm{as $n\to \infty $}.
\end{eqnarray*}
Consequently
\begin{eqnarray*}
e_{n}=\mathcal{O}\big(n^{-1}\big), \ \textrm{as $ n\to \infty$ }
\end{eqnarray*}
 Thus we have shown that 
\begin{eqnarray*}
\big\|dD\mathcal{C}(G)\big\|_{n+1} =\mathcal{O}\big((n-7)\,!\pi^{-n}\big) \ \textrm{as $ n\to \infty $ }.
\end{eqnarray*}

Let $E:=\mathcal{C}(G)$ . Like $G$, the polynomial series $E$ is odd in 
$d$ and the coefficients are odd in $u$. We partition them 
\begin{eqnarray}
E_{n}&=&\alpha_{n}(n-1)\,!\Big(\frac{i}{\pi}\Big)^{n-1}\tau_{n}(u)+\beta_{n-2}(n-3)\,!\Big(\frac{i}{\pi}\Big)^{n-3}\tau_{n-2}(u)\nonumber\\ & &
 {}+\gamma_{n-4}(n-5)\,!\Big(\frac{i}{\pi}\Big)^{n-5}\tau_{n-4}(u)+\overline{E}_{n-6},
\label{5.29}\end{eqnarray}
for odd $n\geq 7$, where $\alpha_{n}$ and $\beta_{n}$ are real
number and also $\overline{E}$ have at most degree $n$ for all $n$.
For the whole series $E$ this is equivalent to 
\begin{eqnarray}
\mathcal{C}(G)=E=E_{1}+d^{2}E_{2}+d^{4}E_{4}+d^{6}\overline{E},
\label{5.30}\end{eqnarray}
where 
\begin{eqnarray*}
E_{1} &=&\sum^{+\infty}_{n=7}\alpha_{n}(n-1)\,!\Big(\frac{i}{\pi}\Big)^{n-1}\tau_{n}(u)d^{n}, \\ 
E_{2}&=&\sum^{+\infty}_{n=5}\beta_{n}(n-1)\,!\Big(\frac{i}{\pi}\Big)^{n-1}\tau_{n}(u)d^{n},
\\ 
E_{3}&=&\sum^{+\infty}_{n=3}\gamma_{n}(n-1)\,!\Big(\frac{i}{\pi}\Big)^{n-1}\tau_{n}(u)d^{n},
\\ 
\overline{E} &=&\sum^{+\infty}_{n=1}\overline{E}_{n}(u)d^{n}.
\end{eqnarray*}
 Theorem \ref{t5.1} and the definition of the norms yields
\begin{eqnarray*}
\alpha_{n}=\mathcal{O}(n^{-7}),
\ \beta_{n}=\mathcal{O}(n^{-5}),\ \gamma_{n}=\mathcal{O}(n^{-3}) 
\ \textrm{and} \ \|D\overline{E}_{n}\|_{n+1}=\mathcal{O}\big((n-1)\,!\pi^{-n} \big).
\end{eqnarray*}
We apply $ \mathcal{C}^{-1}$ to (\ref{5.29}) and  obtain

\begin{equation}
G=\mathcal{C}^{-1}(E_{1})+d^{2}\mathcal{C}^{-1}(E_{2})+d^{4}\mathcal{C}^{-1}(E_{4})+d^{6}\mathcal{C}^{-1}(\overline{E}).
\end{equation}
To the first two summands, theorem \ref{t4.4} applies and yields
\begin{eqnarray*}
\Big\| \{\mathcal{C}^{-1}(E_{1})\}_{n}-\alpha
(n-1)\,!\Big(\frac{i}{\pi}\Big)^{n-1}\tau_{n}(u)\Big\|_{n}&=&\mathcal{O}\big((n-7)\,!\pi^{-n}
\big),\\
\Big\| \{\mathcal{C}^{-1}(E_{2})\}_{n}-\beta
(n-1)\,!\Big(\frac{i}{\pi}\Big)^{n-1}\tau_{n}(u)\Big\|_{n}&=&\mathcal{O}\big((n-5)\,!\pi^{-n}
\big),\\
\Big\| \{\mathcal{C}^{-1}(E_{3})\}_{n}-\gamma
(n-1)\,!\Big(\frac{i}{\pi}\Big)^{n-1}\tau_{n}(u)\Big\|_{n}&=&\mathcal{O}\big((n-3)\,!\pi^{-n}
\big),
\end{eqnarray*}
where 
\begin{eqnarray*}
\alpha=\frac{1}{\pi}\sum^{\infty}_{n=7}\alpha_{n}, \
\beta=\frac{1}{\pi}\sum^{\infty}_{n=5}\beta_{n},\ \gamma=\frac{1}{\pi}\sum^{\infty}_{n=3}\gamma_{n}. 
\end{eqnarray*}
To the last part of (\ref{5.29}), we apply  theorem \ref{t4.3} and obtain
\begin{eqnarray*}
\Big\| \{\mathcal{C}^{-1}(\overline{E})\}_{n}\Big\|_{n}&=&\mathcal{O}\big((n-1)\,!\pi^{-n}\log(n)
\big).
\end{eqnarray*}
Using $G=\mathcal{C}^{-1}(E)$, we have shown 
\begin{eqnarray*}
\bigg\| G_{n}-\alpha
(n-1)\,!\Big(\frac{i}{\pi}\Big)^{n-1}\tau_{n}&-&\beta
(n-3)\,!\Big(\frac{i}{\pi}\Big)^{n-3}\tau_{n-2}-\gamma
(n-5)\,!\Big(\frac{i}{\pi}\Big)^{n-5}\tau_{n-4}\bigg\|_{n}\\
& & {}= \mathcal{O}\Big((n-7)\,!\pi^{-n}\log(n)\Big).
\end{eqnarray*}
If we use (\ref{5.3}) and the relation $\tau_{n+1}=\frac{1}{n}D\tau_{n}$ 
and we apply  theorem \ref{t3.3}, (5), we obtain our final result of this section
\begin{eqnarray*}
\bigg\| F_{n}-\alpha
n\,!\Big(\frac{i}{\pi}\Big)^{n+1}\tau^{\prime}_{n+1}&-&\beta
(n-2)\,!\Big(\frac{i}{\pi}\Big)^{n-1}\tau^{\prime}_{n-1}-\gamma
(n-4)\,!\Big(\frac{i}{\pi}\Big)^{n-3}\tau^{\prime}_{n-3}\bigg\|_{n}= \\
& & {} \mathcal{O}\Big((n-5)\,!\pi^{-n}\log(n+2)\Big).
\end{eqnarray*}
and 
\begin{eqnarray*}
\bigg\| F_{n}-\alpha
(n+1)\,!\Big(\frac{i}{\pi}\Big)^{n+1}\dfrac{\tau_{n+2}}{\tau_{2}}&-&\beta
(n-1)\,!\Big(\frac{i}{\pi}\Big)^{n-1}\dfrac{\tau_{n}}{\tau_{2}}-\gamma
(n-3)\,!\Big(\frac{i}{\pi}\Big)^{n-3}\dfrac{\tau_{n-2}}{\tau_{2}}\bigg\|_{n}\\
&=  & {} \mathcal{O}\Big((n-5)\,!\pi^{-n}\log(n+2)\Big).
\end{eqnarray*}

\section{ Preparation of the functions to construct a quasi-solutions}
In the previous section , we have shown that equation (\ref{2.4}) has a formal solution and
we found an asymptotic approximation of the coefficients of this formal
solution. We will use this to construct the quasi-solution. To 
that purpose, we define the functions 
\begin{eqnarray}
H_{n}(u):=(n+1)\,!\Big(\frac{i}{\pi}\Big)^{n+1}\tau_{n+2}(u)
\label{6.1}\end{eqnarray} 
\begin{eqnarray}
h(t,u):=\sum^{\infty}_{\substack{n=7\\ n\ odd}}H_{n}(u)\frac{t^{n-1}}{(n-1)\,!}.
\label{6.2}\end{eqnarray} 
We rewrite
\begin{eqnarray}
h(t,u)=\Big(\frac{i}{\pi}\Big)^{2}\sum^{\infty}_{\substack{n=6\\ n\ even}}(n+1)(n+2)\Big(\frac{it}{\pi}\Big)^{n}\tau_{n+3}(u).
\label{6.3}\end{eqnarray}
Using proposition (\ref{p3.1})-(4), we obtain
\begin{eqnarray*}
h(t,u)=\frac{-1}{\pi^{2}}\sum^{\infty}_{\substack{n=6\\ n\ even}}\frac{1}{n\,!}\Big(\frac{it}{\pi}\Big)^{n}\frac{d^{n}}{d^{n}\xi}\Big(\frac{d^{2}}{d^{2}\xi}\big(\tanh(\xi)\big)\Big),
\end{eqnarray*}
or equivalently
\begin{eqnarray*}
h(t,u)=\frac{-1}{\pi^{2}}\sum^{\infty}_{\substack{n=6\\ n\ even}}\frac{1}{n\,!}\Big(\frac{it}{\pi}\Big)^{n}\frac{d^{n}}{d^{n}\xi}\Big(g(\xi)\Big),
\end{eqnarray*}
where \ $g(\xi)=-2\big(\tanh(\xi)-\tanh(\xi)^{3}\big)$ and $\xi=\xi(u)=artanh(u).$\\
and hence that 
\begin{eqnarray*}
h(t,u)=\frac{-1}{\pi^{2}}\sum^{\infty}_{\substack{n=0\\ n\ even}}\frac{1}{n\,!}\Big(\frac{it}{\pi}\Big)^{n}\frac{d^{n}}{d^{n}\xi}\Big(g(\xi)\Big)+\mu(t,u),
\end{eqnarray*}
where $$ \mu(t,u):=-\frac{2\, u}{\pi^{2}}\tau_2(u)+\frac{4}{\pi^4}(3u^3-2u)\tau_2(u)\,t^2-\frac{2}{3\pi^6}(45u^5-60u^3+17u)\tau_2(u)\,t^4.$$
Thus we can  write
\begin{eqnarray}
h(t,u)=\frac{-1}{2\pi^{2}}\Big[g\big(\xi+\frac{it}{\pi}\big)+g\big(\xi-\frac{it}{\pi}\big)\Big]+ \mu(t,u).
\label{6.4}\end{eqnarray}
This gives
\begin{eqnarray}
h(t,u)&=&(1-u^2)\Bigg[\frac{(4u^{3}-4u)\sin\big(\frac{t}{\pi}\big)^{4}-(6u^{3}-2u)\sin\big(\frac{t}{\pi}\big)^{2}+2u}{\pi^{2}\Big(\cos\big(\frac{t}{\pi}\big)^{2}+u^{2}\sin\big(\frac{t}{\pi}\big)^{2}\Big)^{3}}\Bigg]\nonumber\\
&+& \mu(t,u). 
\label{6.5}\end{eqnarray}
For fixed real $u$ the function $h(.,u) $ is analytic in $|t|<\rho$, where $\rho=\frac{\pi^2}{2}$. In the subsequent definition, 
we consider real values  of $u$,   $0< u\leq 1$, here
$h(.,u)$ is also analytic with respect to $t$ on the positive real axis.

Now we define the function $\mathcal{H}(d,u)$ by
\begin{eqnarray}
\mathcal{H}(d,u):= \int^{+\infty}_{0}e^{-\frac{t}{d}}h(t,u) d\,t, \quad \text{for}  \ \ \ (0< u\leq 1).
\label{6.6}\end{eqnarray}

The function $\mathcal{H}(d,.)$ is real analytic; they can be continued analytically
to the interval $-1<u\leq1$ in the following way. Choose some positive number $M$ and
let $\Gamma_1$ the path consisting of the segment from $0$ to $Mi$ and
of the ray $t\mapsto t+Mi, t\geq0$. Let $\Gamma_2$ the symmetric path that could 
also be obtained using $-M$ instead of $M$. Recalling (\ref{6.4}), we can also define
\begin{eqnarray*}
\mathcal{H}(d,u):&=&-\frac{1}{\pi^2}
\left[\int_{\Gamma_2}e^{-\frac{t}{d}}
g\big(\xi+\frac{it}{\pi}\big)\,dt+
\int_{\Gamma_1}e^{-\frac{t}{d}}g\big(\xi-\frac{it}{\pi}\big)\,dt
\right]\nonumber\\
&+&\mu_{1}(d,u), 
\end{eqnarray*}
where $$\mu_{1}(d,u):=\int_{0}^{\infty}e^{-\frac{t}{d}}\mu(t,u)dt,$$ 
for $-\tanh(\frac2\pi M)<u\leq1$, where $\xi=\mbox{artanh}(u)$, 
because the singularities of $\tanh$ are $i(\frac\pi2+n\pi)$, 
$n$ integer. As $M$ is arbitrary, this defines the analytic continuation
of $\mathcal{H}(d,.)$ for $-1<u\leq1$.

 In the sequel we consider $u_0\in]-1,0]$.
\begin{lem}  If we consider the above function $\mathcal{H}(d,u)$ and the operators defined in (\ref{2.22}). 
Then, for $u_0<u\leq1$
\begin{itemize}
\item 1.  $$\mathcal{S}_{2}(\mathcal{H})=(1-u^2) \mathcal{D}_{1}(d,u),$$
 where  the function $ \mathcal{D}_{1}(d,u)$ is analytic,  beginnings with $d^{8}$.
\item2. 
\begin{eqnarray*}
\mathcal{C}_{2}(\mathcal{H})=-\mathcal{H}(d,u)+(1-u^2) \mathcal{D}_{2}(d,u),
\end{eqnarray*}
  where  the function $ \mathcal{D}_{2}(d,u)$ is analytic,  beginnings with $d^{7}$.
\item 3 
\begin{eqnarray*}
\mathcal{S}_{2}\Big(\dfrac{1}{\tau_2(u)}\mathcal{H}\Big)=- \dfrac{2u\varepsilon}{\tau_2(u)}\mathcal{H}(d,u) +\mathcal{D}_3(d,u)
\end{eqnarray*}
where the function $\mathcal{D}_3(d,u)$ is analytic beginnings with $d^8$
\end{itemize}
\label{lem6.1.1}\end{lem}
\noindent{\bf Proof.}
(1)-
If we replace $u$ by $T^{+}$  and $T^{-}$ in (\ref{6.5})  we obtain
\begin{eqnarray*}
\mathcal{H}(d,T^{+})-\mathcal{H}(d,T^{-})= \int^{+\infty}_{0}e^{-\frac{t}{d}}\Big(h(t,T^{+})-h(t,T^{-})\Big) d\,t,
\end{eqnarray*}
this implies
\begin{eqnarray}
\mathcal{H}(d,T^{+})-\mathcal{H}(d,T^{-})=-\frac{1}{2\pi^{2}}\big(\mathcal{I}^{+}-\mathcal{I}^{-}\big)+\sigma(d,T^{+})-\sigma(d,T^{-}),
\label{6.7}\end{eqnarray}
where 
\begin{eqnarray*}
\mathcal{I}^{+}&=&\int^{+\infty}_{0}e^{-\frac{t}{d}}\Big(g\big(\xi(T^{+})+\frac{it}{\pi}\big)+g\big(\xi(T^{+})-\frac{it}{\pi}\big)\Big)d\,t, 
 \\
\mathcal{I}^{-}&=&\int^{+\infty}_{0}e^{-\frac{t}{d}}\Big(g\big(\xi(T^{-})+\frac{it}{\pi}\big)+g\big(\xi(T^{-})-\frac{it}{\pi}\big)\Big)d\,t,\\
 \sigma(d,u) &=&\int^{+\infty}_{0}e^{-\frac{t}{d}}\mu(t,u) d\,t.
\end{eqnarray*}
Using \, $\xi(T^{\pm})=\xi\pm d$, where $\xi$ and $u$ are  coupled 
by $\xi=\text{artanh}(u)$, we obtain
 \begin{eqnarray*}
\mathcal{I}^{+}=\int^{+\infty}_{0}e^{-\frac{t}{d}}g\big(\xi+d+\frac{it}{\pi}\big)d\,t +\int^{+\infty}_{0}e^{-\frac{t}{d}}g\big(\xi+d-\frac{it}{\pi}\big)d\,t.
\end{eqnarray*}
If we  substitute $t+\pi i\,d $ in the first
part,\ $t-\pi i\,d $ in the second part, we obtain 
 \begin{eqnarray*}
\mathcal{I}^{+}=-\int^{+\infty+\pi id}_{-\pi id}e^{-\frac{t}{d}}g\big(\xi+\frac{it}{\pi}\big)d\,t -\int^{+\infty+\pi id}_{\pi id}e^{-\frac{t}{d}}g\big(\xi-\frac{it}{\pi}\big)d\,t.
\end{eqnarray*}
We apply Cauchy's theorem 
\begin{eqnarray*}
\mathcal{I}^{+}=-\int^{+\infty}_{0}e^{-\frac{t}{d}}\Big(g\big(\xi+\frac{it}{\pi}\big)+g\big(\xi-\frac{it}{\pi}\big)\Big)d\,t
&+&\int^{-\pi id}_{0}e^{-\frac{t}{d}}g\big(\xi+\frac{it}{\pi}\big)d\,t\\
 &+& \int^{\pi id}_{0}e^{-\frac{t}{d}}g\big(\xi-\frac{it}{\pi}\big)d\,t.
\end{eqnarray*}
 Substitute $t=-i\,s\,d $ in the second
part,\ $t= i\,s\,d $ in the third part, we obtain 
\begin{eqnarray*}
\mathcal{I}^{+}=-\int^{+\infty}_{0}e^{-\frac{t}{d}}\Big(g\big(\xi+\frac{it}{\pi}\big)+g\big(\xi-\frac{it}{\pi}\big)\Big)d\,t+2d\int^{\pi}_{0}\sin(s) g\big(\xi+\frac{s\,d}{\pi}\big)d\,s.
\end{eqnarray*}
This implies 
\begin{eqnarray}
\mathcal{I}^{+}&=&-\int^{+\infty}_{0}e^{-\frac{t}{d}}\Big(g\big(\xi+\frac{it}{\pi}\big)+g\big(\xi-\frac{it}{\pi}\big)\Big)d\,t-4d(1-u^{2})\mathcal{D}^{+}(d,u),
\label{6.8}\end{eqnarray} 
where 
\begin{eqnarray*}
\mathcal{D}^{+}(d,u)=\int^{\pi}_{0}\sin(s)\frac{\big(u+\tanh(\frac{s\,d}{\pi})\big)\big(1-\tanh(\frac{s\,d}{\pi})^{2}\big)}{\big(1+u\tanh(\frac{ds}{\pi})\big)^{3}}d\,s.
\end{eqnarray*} 
We can also use the same method for $\mathcal{I}^{-}$ and obtain
\begin{eqnarray}
\mathcal{I}^{-}&=&-\int^{+\infty}_{0}e^{-\frac{t}{d}}\Big(g\big(\xi+\frac{it}{\pi}\big)+g\big(\xi-\frac{it}{\pi}\big)\Big)d\,t-4d(1-u^{2})\mathcal{D}^{-}(d,u),
\label{6.9}\end{eqnarray} 
where 
\begin{eqnarray*}
\mathcal{D}^{-}(d,u)=\int^{\pi}_{0}\sin(s)\frac{\big(u-\tanh(\frac{ds}{\pi})\big)\big(1-\tanh(\frac{s\,d}{\pi})^{2}\big)}{\big(1-u\tanh(\frac{s\,d}{\pi})\big)^{3}}d\,s.
\end{eqnarray*} 
Consequently
\begin{eqnarray*}
\mathcal{H}(d,T^{+})-\mathcal{H}(d,T^{-})&=&\frac{2d(1-u^{2})}{\pi^{2}}\Big(\mathcal{D}^{+}(d,u)-\mathcal{D}^{-}(d,u)\Big)+\sigma(d,T^{+})-\sigma(d,T^{-}).
\end{eqnarray*}
Using (\ref{2.22}), we find
\begin{eqnarray}
\mathcal{S}_{2}(\mathcal{H})=(1-u^2)\,\mathcal{D}_{1}(d,u),
\label{6.10}\end{eqnarray}
where 
$\mathcal{D}_{1}(d,u)$ is  analytic, beginnings with $d^8.$

(2). Using the same method, we obtain
\begin{eqnarray*}
\mathcal{H}(d,T^{+})+\mathcal{H}(d,T^{-})=-2\mathcal{H}(d,u)+2(1-u^{2})\mathcal{D}_{2}(d,u).
\end{eqnarray*}
  With (\ref{2.22}) we find
\begin{eqnarray*}
\mathcal{C}_{2}(\mathcal{H})=-\mathcal{H}(d,u)+(1-u^{2})\mathcal{D}_{2}(d,u).
\end{eqnarray*}
$\mathcal{D}_{2}(d,u)$ is  analytic, beginnings with $d^7.$

(3). Using (\ref{5.7}), (\ref{2.8}), (\ref{2.22}) and (1), (2) of this   lemma, we obtain
\begin{eqnarray*}
\mathcal{S}_{2}\Big(\dfrac{1}{\tau_2(u)}\mathcal{H}\Big)&=& \mathcal{S}_{2}\Big(\dfrac{1}{\tau_2(u)}\Big)\mathcal{C}_{2}\big(\mathcal{H}\big)+\mathcal{C}_{2}\Big(\dfrac{1}{\tau_2(u)}\Big)\mathcal{S}_{2}\big(\mathcal{H}\big)\\
&=&- \dfrac{2u\varepsilon}{\tau_2(u)}\mathcal{H}(d,u) +\mathcal{D}_3(d,u)
\end{eqnarray*}
where $$\mathcal{D}_3(d,u):=\big(1+(1+u^2)\sinh(d)^2\big)(1-u^2)\mathcal{D}_{1}(d,u)+2u\varepsilon(1-u^2)\mathcal{D}_{2}(d,u).$$
\begin{propo}
We have
\begin{enumerate}
\item For $u_0<u\leq 1$
\begin{equation}
\mathcal{H}(d,u)\sim\sum^{\infty}_{\substack{n=7\\ n\  odd}}(n+1)\,!\Big(\frac{i}{\pi}\Big)^{n+1}\tau_{n+2}(u)d^{n},\ 
\textrm{as} \ d\searrow 0,
\label{6.11}\end{equation} 

\item $\big|\frac{\partial\mathcal{H}}{\partial u}(d,u)\big|\leq Kd$\quad  for \
  $u_0\leq u\leq 1\ $ $(d>0).$
\end{enumerate}
\label{p6.1}\end{propo}
\textbf{Proof.}
\begin{enumerate}
\item To prove (1) we use Watson's lemma and (\ref{6.2}).
  if$\ |t-k\pi^{2}|<\frac{\pi^{2}}{2}$,\quad for  $u_0\leq u< 1$ 
\item 
\begin{itemize}
\item For $u_0\leq u\leq \frac{1}{2}$ 
\begin{eqnarray*}
\Big|\mathcal{H}(d,u)\Big|=\frac{1}{d^{6}}\int_{0}^{\infty}e^{-t/d} h^{(-6)}(t,u)\,dt,
\end{eqnarray*}
where  $h^{(-6)}(t,u)$ satisfies
\begin{equation*}
\Big(\frac{\partial}{\partial
  t}\Big)^{6}h^{(-6)}=\frac{\partial{h}}{\partial u}\quad \text{and} \Big(\frac{\partial}{\partial
  t}\Big)^{k}h^{(-6)}(0,u)=0,\quad \text{for}\, k=0,...,5.
\end{equation*}
 
Because of
\begin{eqnarray*}
h^{(-6)}(t,u)\leq
K\frac{t^{7}}{7\,!} \qquad \textrm{for}\  t\geq 0,
\end{eqnarray*}
we obtain 
\begin{eqnarray*}
\Big|\frac{\partial\mathcal{H}}{\partial
  u}(d,u)\Big|&\leq&\frac{K}{7\,!d^{6}}\int^{\infty}_{0}e^{-\frac{t}{d}}t^{7}dt.
\end{eqnarray*}
This implies 
\begin{eqnarray*}
\Big|\frac{\partial\mathcal{H}}{\partial
  u}(d,u)\Big|
\leq Kd  \ \quad \quad (d>0).
\end{eqnarray*}
\item for $\frac{1}{2}\leq u\leq 1$ 
\begin{eqnarray*}
\Big|\frac{\partial\mathcal{H}}{\partial
  u}(d,u)\Big|&=&\Big|\int^{\infty}_{0}e^{-\frac{t}{d}}\frac{\partial{h}}{\partial
  u}(d,u)dt\Big|, \\
&\leq&K\int^{\infty}_{0}e^{-\frac{t}{d}}dt,\\
&\leq&Kd  \ \quad \quad (d>0).\ \square\
\end{eqnarray*}
\end{itemize}
\end{enumerate}

Let us  consider a sequence $R_{n}(u)$ of polynomials of degree at most $n$,
such that 
\begin{eqnarray*}
\|R_{n}\|_{n}=O\Big((n-5)\,!\pi^{-1}\log(n)\Big)
\end{eqnarray*}
\begin{lem}
Let polynomial series  difine $R\in \mathcal{Q}$ with the above
estimate be given and define
\begin{eqnarray*}
r(t,u):&=&\sum_{n=7}^{\infty}R_{n}(u)\frac{t^{n-1}}{(n-1)\,!}, \  \quad
(t\in \mathbb{C}, |t|\leq\frac{\pi^{2}}{2}, u_0\leq u\leq 1 ),\\
r(t,u):&=&r\big(\frac{\pi^{2}}{2},u\big)+\big(t-\frac{\pi^{2}}{2}\big)\frac{\partial
  r}{\partial t}\big(\frac{\pi^{2}}{2},u\big), \quad
( t>\frac{\pi^{2}}{2}, u_0\leq u\leq 1 ),\\
\mathcal{R}(d,u):&=&\int_{0}^{\infty}e^{-\frac{t}{d}}\,r(t,u)dt. 
\end{eqnarray*}
Then
\begin{enumerate}
\item $r$ is continuously differentiable function on the set $B$ of
all $(t,u)$ such that $u$ satisfies $u_0\leq u\leq 1$ and $t$ is a
complex number and satisfies  $|t|\leq\frac{\pi^{2}}{2}$ or $
t>\frac{\pi^{2}}{2}$. The restriction of $r$ to  $u_0\leq u\leq
1,|t|\leq\frac{\pi^{2}}{2}$ is twice continuously differentiable. For
fixed $u_0\leq u\leq 1$ the function $r(t,u)$ is analytic in
$|t|<\frac{\pi^{2}}{2}.$
\item $\mathcal{R}(d,u)$ is continuous, partially differentiable with
  respect to $u$, has  continuous partial derivative and 
\begin{eqnarray}
\mathcal{R}(d,u)\sim\sum^{\infty}_{n=7}R_{n}(u)d^{n},\ 
\textrm{as} \ d\searrow 0.
\label{6.12}\end{eqnarray}
\item  $\big|\mathcal{R}(d,u)\big|\leq Kd^{3},\
  \big|\frac{\partial\mathcal{R}}{\partial u}(d,u)\big|\leq Kd^{3},$\quad  for \
  $u_0\leq u\leq 1\ $ $(d>0).$
\end{enumerate}
\label{l6.3}\end{lem}
The importance of our definition of $\mathcal{R}$ lies in a certain
compatibility with insertion of the functions $T^{+}, T^{-}$ for
$u$. First let
\begin{eqnarray*}
\sum^{\infty}_{n=7}R_{n}^{+}(u)d^{n}&:=&\sum^{\infty}_{n=7}R_{n}\big(T^{+}(d,u)\big)d^{n},\\
\sum^{\infty}_{n=7}R_{n}^{-}(u)d^{n}&:=&\sum^{\infty}_{n=7}R_{n}\big(T^{-}(d,u)\big)d^{n}.
\end{eqnarray*}
We obtain a new sequences $R^{+}_{n}(u), R^{-}_{n}(u)$ of polynomials
of degree at most $n$.  Theorem \ref{t4.2},
and (\ref{4.2}) imply 
\begin{eqnarray*}
\big\|R^{+}_{n}(u)\big\|_{n}&=&O\Big((n-5)\,!\pi^{-n}\log(n)\Big)\\
\big\|R^{-}_{n}(u)\big\|_{n}&=&O\Big((n-5)\,!\pi^{-n}\log(n)\Big)
\end{eqnarray*}
Therefore we can use lemma \ref{l6.3} for $R^{+}_{n}(u),\  R^{-}_{n}(u)$,
and obtain  functions $\mathcal{R}^{+}(d,u)$, $ \mathcal{R}^{-}(d,u)$.
\begin{theo}
There is a positive constant $K$ independent of $ d,u$ such that
\begin{eqnarray*}
\big|\mathcal{R}^{+}(d,u)-\mathcal{R}(d,T^{+})\big|&\leq&
Kd^{3}e^{-\frac{\pi^{2}}{2d}}\ \ \  \text{for}\ \ \ (0<d<d_0,  u_0<u<1) \\
\big|\mathcal{R}^{-}(d,u)-\mathcal{R}(d,T^{-})\big|&\leq& Kd^{3}e^{-\frac{\pi^{2}}{2d}}\ \ \  \text{for}\ \ \ (0<d<d_0,  u_0<u<1) .
\end{eqnarray*}
\label{t6.4}\end{theo}
\textbf{Proof}
The proof is exactly the one of \cite{SV}.
\begin{defi}
Let $\mathcal{D}(d,u)$ be a function defined for $0<d<d_{0}$ and
$u_0<u<1$. We say that $\mathcal{D}(d,u)$ has property \textmd{G} if 
\begin{eqnarray*}
\mathcal{D}(d,u)=\int_{0}^{\infty}e^{-\frac{t}{d}}\,q(t,u)dt \quad
(0<d<d_{0},\ u_0<u<1 )
\end{eqnarray*}
is the Laplace transform of some function $q(t,u)$ that has the
following proprties:  
\begin{enumerate}
\item  $q(t,u)$ is defined if $u_0<u<1$  and either $t$ is complex and
  $|t|<\frac{\pi^{2}}{2}$ or $t$ is real and $t\geq 0$,
\item $q(t,u)$ is analytic in $|t|<\frac{\pi^{2}}{2}$ for $ u_0<u<1$,
\item $q(t,u)$ restricted to $0\leq t<\frac{\pi^{2}}{2}$ or
  $t\geq\frac{\pi^{2}}{2}$ is continuous and the
  $\lim_{t\to\frac{\pi^{2}}{2} } q(t,u)$ exists for every  $ u_0<u<1$,
\item there is a positive constant $K$  such that
\begin{eqnarray*}
 |q(t,u)|\leq K e^{Kt},\ \  \text{for}\ \ t\geq 0,\  u_0<u<1, \,\,
(0<d<d_{0},\ u_0<u<1 )
\end{eqnarray*}
\end{enumerate}
The set of all function with the property \textmd{G} will be denote by
$\mathcal{G}$
\label{d6.5}\end{defi}
\begin{lem} For $u_0<u\leq 1$
\begin{enumerate}
\item If $\mathcal{H}(d,u)$ is the function defined in (\ref{6.6}) then 
$$d^4\mathcal{H}(d,u)=(1-u^2)\tilde{\mathcal{H}}(d,u)+\mathcal{O}\Big((1-u^2)e^{-\frac{\pi^2}{2d}}\Big)$$
where $\tilde{\mathcal{H}}(d,u)$ has property \textmd{G}

\item  Let $k$ be a psitive integer. If  $\mathcal{D}_{1}, \mathcal{D}_{2}$ have  property \textmd{G} and  their first terms in the Taylor development at $d = 0$, begin with $d^k$  then 
$$\mathcal{D}_{1}(d,u)\mathcal{D}_{2}(d,u)=d^k\mathcal{D}(d,u)+\mathcal{O}\Big(d^ke^{-\frac{\pi^2}{2d}}\Big),$$ where $\mathcal{D}(d,u)$ has  property \textmd{G}
\item If $\mathcal{H}(d,u)$ is the function defined in (\ref{6.6}) then
$$\mathcal{H}(d,u)^2=(1-u^2)^2 d^3\mathcal{E}(d,u)+\mathcal{O}\Big((1-u^2)^2 d^3e^{-\frac{\pi^2}{2d}}\Big)$$
where $\mathcal{E}(d,u)$ has property \textmd{G}
\item Any function  $\mathcal{D}(d,u)$ analytic in a neighborhood of $d=0$ has  property \textmd{G} if  $\mathcal{D}(0,u)=0$ for all $u$,
\item If $\mathcal{R}(d,u)$ is defined by lemma \ref{l6.3} then 
  $\dfrac{1}{d^{2}}\mathcal{R}(d,u)$  has  property \textmd{G}
\item If  $\mathcal{D}_{1}, \mathcal{D}_{2} $ have  property
  \textmd{G} then so do $\mathcal{D}_{1}+\mathcal{D}_{2},
 \  \mathcal{D}_{1}-\mathcal{D}_{2}$ and\, $\mathcal{D}_{1}\cdot\mathcal{D}_{2} $
\item If  $\mathcal{D}(d,u)$ has  property \textmd{G} then 
\begin{eqnarray}
\big|\mathcal{D}(d,u)\big|\leq Kd  \quad \big(0<d<\frac{1}{K})
\label{6.13}\end{eqnarray}
with some constant $K>0$ independent of $u$.
\end{enumerate}
\label{l6.6}\end{lem}
\noindent{\bf Proof.}

\begin{enumerate}
\item   For $u>0$,
 we have 
\begin{eqnarray}
d^4 \mathcal{H}=(1-u^2) \int_{0}^{\infty} e^{-\frac{t}{d}}g_4(t,u) dt
\label{fp1}\end{eqnarray}
where $$g_4(t,u)=\dfrac{1}{(1-u^2)}\int_{0}^{t}\int_{0}^{\theta}\int_{0}^{\upsilon}\int_{0}^{\tau}h(s,u) ds\, d\tau\, d\upsilon\ d\theta $$
$g_4(t,u)$ has a logarithmic singularity at $t_k(s)=(2k+1)\dfrac{\pi^2}{2}\pm d\frac{\pi\, s }{\varepsilon}i$ for $(k\geq 0, s>0)$. it is analytic in $|t|<\frac{\pi^2}{2}$ and $\lim_{t\to\frac{\pi^2}{2}}g_4(t,u)$ exists. 

If  we put 
$$\tilde{\mathcal{H}}(d,u)=\int_{0}^{\infty} e^{-\frac{t}{d}}\tilde{g}_4(t,u) d\,t $$
where $$\tilde{g}_4(t,u)=
\begin{cases}
g_4(t,u),          &\text{if $t\leq \frac{\pi^2}{2}$}\\
g_4\big(\frac{\pi^2}{2},u\big),          &\text{if $t\geq \frac{\pi^2}{2}$}
\end{cases}
$$
then $\tilde{\mathcal{H}}(d,u)$ has property \textmd{G} and 
$$d^4\mathcal{H}(d,u)=(1-u^2)\tilde{\mathcal{H}}(d,u)+\mathcal{O}\Big((1-u^2)e^{-\frac{\pi^2}{2d}}\Big). $$

(ii)- For $-u_0 <u<1$, where $0<u_0 <1$,  we have 
\begin{eqnarray*}
\mathcal{H}(d,u)&=& \int_{0}^{\infty e^{i\varphi}} e^{-\frac{t}{d}}h(t,u) d\,t+2\pi i\sum_{k\geq 0}Res\Big(e^{-\frac{t}{d}}h(t,u),t_k(s)\Big)\\
&=&\int_{0}^{\infty e^{i\varphi}} e^{-\frac{t}{d}}h(t,u) d\,t+\mathcal{O}\bigg(\frac{1}{d^2}(1-u^2)e^{-\frac{\pi^2}{2d}}\bigg)
\end{eqnarray*}
where $\dfrac{\pi}{2}<\varphi<\dfrac{\pi}{4}$. For $0<u<u_0$, this formula coincides with the formula $$\int_{0}^{\infty} e^{-\frac{t}{d}}h(t,u) d\,t $$
and extends it by real analytic continuation  for $-u_0<u<0.$

This implis 
\begin{eqnarray*}
d^4\mathcal{H}(d,u)&=&(1-u^2)\int_{0}^{\infty e^{i\varphi}} e^{-\frac{t}{d}}g_4(t,u) d\,t+\mathcal{O}\Big((1-u^2)d^2e^{-\frac{\pi^2}{2d}}\Big)
\end{eqnarray*}
we obtain 
\begin{eqnarray*}
d^4\mathcal{H}(d,u)&-&(1-u^2)\tilde{\mathcal{H}}(d,u)=(1-u^2)\int_{\Gamma} e^{-\frac{t}{d}}g_2(t,u) d\,\Gamma\\
&+&\mathcal{O}\Big((1-u^2)e^{-\frac{\pi^2}{2d}}\Big),
\end{eqnarray*}
where $\Gamma$ is the path following the real line from infinity to $(\frac{\pi^2}{2},0)$, then along the vertical line from $(\frac{\pi^2}{2},0)$ to $\Big(\frac{\pi^2}{2},\frac{\pi^2}{2}\tan(\varphi)\Big)$ and finally along the line $y=\tan(\varphi) x$ from  $\Big(\frac{\pi^2}{2},\frac{\pi^2}{2}\tan(\varphi)\Big)$ to infinity.

Since $g_4(.,u)$ is bounded on $\Gamma$, then 
$$\Big|d^4\mathcal{H}(d,u)-(1-u^2)\tilde{\mathcal{H}}(d,u)\Big|\leq K (1-u^2)e^{-\frac{\pi^2}{2d}},$$
where $K$ is positive constant. Finally 
$$d^4\mathcal{H}(d,u)=(1-u^2)\tilde{\mathcal{H}}(d,u)+\mathcal{O}\Big((1-u^2)e^{-\frac{\pi^2}{2d}}\Big) \ \ \ \text{for} \ \ \  (0<u_0<u\leq1)$$



\item We assume that  $\mathcal{D}_{1}(d,u), \mathcal{D}_{2}(d,u)$ have  property \textmd{G} and their first terms in the Taylor development at $d=0$ begin with $d^k$. Then 
\begin{eqnarray*}
\mathcal{D}_{1}= \int_{0}^{\infty} e^{-\frac{t}{d}}f(t,u) d\,t\\
\mathcal{D}_{2}= \int_{0}^{\infty} e^{-\frac{t}{d}}g(t,u) d\,t
\end{eqnarray*}
where $f(t,u), g(t,u)$ are analytic in $|t|<\frac{\pi^2}{2}$ and $f(t,u)=O(t^{k-1}), g(t,u)=O(t^{k-1})$. 
\begin{eqnarray}
\mathcal{D}_{1}(d,u)\mathcal{D}_{2}(d,u)= \int_{0}^{\infty} e^{-\frac{t}{d}}(f*g)(t,u) d\,t
\end{eqnarray}
Since 
\begin{eqnarray*}
 h(t,u)&= &(f*g)(t,u) =\int_{0}^{s}f(t,u)g(t-s,u)ds \\
&=&\int_{0}^{t}f(t-s,u)g(s,u)ds\\
&=& \int_{t/2}^{t}f(s,u)g(t-s,u)ds +\int_{t/2}^{t}f(t-s,u)g(s,u)ds.
\end{eqnarray*}
 For $t<\frac{\pi^2}{2}$, the function $h(t,u)$ is  $k$ times differentiable  with respect to $t$  and 
\begin{eqnarray*}
 h'(t,u)&= & f(t,u)g(0,u)+f(0,u)g(t,u)-f(\frac{t}{2},u)g(\frac{t}{2},u)\\
&+& \int_{t/2}^{t}f(s,u)g'(t-s,u)ds +\int_{t/2}^{t}f'(t-s,u)g(s,u)ds\\
&=&\int_{t/2}^{t}f(s)g'(t-s)ds +\int_{t/2}^{t}f'(t-s,u)g(s,u)ds\\
&-&f(\frac{t}{2},u)g(\frac{t}{2},u)\\
 h^{(k)}(t,u)&= &\int_{t/2}^{t}f(s)g^{(k-1)}(t-s)ds +\int_{t/2}^{t}f^{(k-1)}(t-s,u)g(s,u)ds\\
&- &\sum_{n=0}^{^k-1}f^{(n)}(\frac{t}{2},u)g^{(k-1-n)}(\frac{t}{2},u)
\end{eqnarray*}
Observe that $ h^{(k)}(t,u)$ is continuous on $[0,\pi^2[$, it is analytic for $|t|<\frac{\pi^2}{2}$. If we put 
 $$\tilde{h}(t,u)=
\begin{cases}
h(t,u),          &\text{if $t< \frac{\pi^2}{2}$}\\
\Big(t-\frac{\pi^2}{2}\Big)^k.         &\text{if $t\geq \frac{\pi^2}{2}$}
\end{cases}
$$
then 
$$ \int_{0}^{\infty} e^{-\frac{t}{d}}\tilde{h}^{(k)}(t,u) d\,t,$$  has property \textmd{G} and 
\begin{eqnarray*}
 (\mathcal{D}_{1}\cdot\mathcal{D}_{2})(d,u)&=&  \int_{0}^{\infty} e^{-\frac{t}{d}}h(t,u) d\,t=\int_{0}^{\infty} e^{-\frac{t}{d}}\tilde{h}(t,u) d\,t+\mathcal{O}\Big(d^ke^{-\frac{\pi^2}{2d}}\Big) \\
&=&d^k\int_{0}^{\infty} e^{-\frac{t}{d}}\tilde{h}^{(k)}(t,u) d\,t+\mathcal{O}\Big(d^ke^{-\frac{\pi^2}{2d}}\Big)\\
&=&d^k\mathcal{D}(d,u)+\mathcal{O}\Big(d^ke^{-\frac{\pi^2}{2d}}\Big)
\end{eqnarray*}
where $\mathcal{D}(d,u)$ has property \textmd{G}.
\item We have 
$$\mathcal{H}(d,u)^2=\frac{1}{d^8}(d^4\mathcal{H})\cdot ( d^4\mathcal{H})$$ 
Using (1) and (2) of this lemma, we obtain immediately  the result.
\end{enumerate}
 For $0\leq u\leq 1,$ the proof of   (4), (5), (6)  and (7) is exactly the one of  \cite{SV}. This proof is valid for  $u_0<u\leq 1.$ 


\section{Approximate solution $ \mathcal{A}(d,u)$ of the functional
  equation (\ref{2.4})}
In this section we use the formal solution of  (\ref{2.4}). The estimates of 
section 5 for its coefficients and the preliminary results of section
6 allow to construct the  function  $ \mathcal{A}(d,u)$ that almost satisfies (\ref{2.4}),
except an error which is exponentially small as $d\to 0$.

In theorem \ref{t2.1} we found that (\ref{2.4}) has a uniquely determined formal
power series solution 
\begin{eqnarray}
A(d,u)=\frac{1}{\varepsilon}B(d,u)=\frac{1}{\varepsilon}\bigg(U+\sum_{\substack{n=7\\ 
    n\, \textrm{odd}}}^{\infty}F_{n}(u)d^{n}\bigg),
\label{7.1}\end{eqnarray}
where 
\begin{eqnarray*}
U(d,u)&=&\varepsilon u
+(u-u^{3})d^{3}+(\frac{10}{3}u^{5}-\frac{16}{3}u^{3}+2u)d^{5},
\\
F_{n}&=&\alpha n\,!\Big(\frac{i}{\pi}\Big)^{n+1}\tau^{\prime}_{n+1}(u)+\beta(n-2)\,!\Big(\frac{i}{\pi}\Big)^{n-1}\tau^{\prime}_{n-1}(u)\\&+&\gamma (n-4)\,!\Big(\frac{i}{\pi}\Big)^{n-4}\tau^{\prime}_{n-3}(u)+R_{n}(u),\\
\big\|R_{n}\big\|_{n}&=&O\Big((n-5)\,!\pi^{-1}\log(n)\Big), \quad
(n\geq 7).\\
\end{eqnarray*}
Consequently
\begin{equation*}
B(d,u)= \big(\alpha+\beta d^{2}+ \gamma d^{4}\big)\sum_{\substack{n=7\\  n
    \,\textrm{odd}}}^{\infty}(n+1)\,!\Big(\frac{i}{\pi}\Big)^{n+1}\dfrac{\tau_{n+2}(u)}{\tau_{2}(u)}d^{n} +\sum_{\substack{n=7\\  n \,\textrm{odd}}}^{\infty}R_{n}(u)d^{n} +U(d,u),
\end{equation*}
Now we set 
\begin{equation}
\mathcal{B}(d,u)= \big(\alpha+\beta d^{2}+ \gamma
d^{4}\big)\dfrac{1}{\tau_{2}(u)} \mathcal{H}(d,u)+\mathcal{R}(d,u) +U(d,u),\  (d>0, u_0<u<1),
\label{7.2}\end{equation}
 where $\mathcal{H}(d,u)$ is defined in (\ref{6.6}) and
$\mathcal{\tilde{R}}(d,u)$ is the function corresponding to
$R_{n},\  n=7,9,..$ By lemma \ref{l6.3} and proposition \ref{p6.1}, we have
\begin{equation}
\mathcal{B}(d,u)\sim B(d,u),\quad \textrm {as} \ d\searrow 0\  \textrm
{for every}\ \ \   u_0<u<1.
\label{7.3}\end{equation}
\begin{theo}
We have 
\begin{equation*}
\Big|\mathcal{A}(d,T^{+})-\mathcal{A}(d,T^{-})-f\big(\varepsilon,\mathcal{A}(d,u)\big)\Big|\leq 
\frac{K}{d}e^{-\frac{\pi^{2}}{2d}},\quad \textrm{for}\ (0<d<d_{0}, 
u_0\leq u<1),
\end{equation*}
where $K$ is a constant independent of $d, u$ and 
\begin{equation}
\mathcal{A}(d,u)=\frac{1}{\varepsilon}\mathcal{B}(d,u)
\label{7.4}\end{equation}
\label{t7.1}\end{theo}

\noindent\textbf{Proof.}
We set 
\begin{equation}
\mathcal{T}(d,u)=\mathcal{S}_{2}(\mathcal{B})-\varepsilon^{2}+\mathcal{B}(d,u)^{2}.
\label{7.5}\end{equation}
then 
\begin{eqnarray*}
\mathcal{T}(d,u)&=&\big(\alpha+\beta d^{2}+ \gamma
d^{4}\big)\mathcal{S}_{2}\Big(\dfrac{1}{\tau_{2}(u)} \mathcal{H}\Big)-\varepsilon^{2}+\mathcal{B}(d,u)^{2}\\
&+&\frac{1}{2}\mathcal{R}(d,T^+) -\frac{1}{2}\mathcal{R}(d,T^-) +\mathcal{S}_{2}(U)
\end{eqnarray*}
 With (\ref{7.2}) and  lemma \ref{lem6.1.1}, this implies 
\begin{equation*}
\mathcal{T}(d,u)=\mathcal{T}_{0}(d,u)\cdot\mathcal{H}(d,u)+\mathcal{T}_{1}(d,u)+\mathcal{T}_{2}(d,u),
\end{equation*}
where
\begin{eqnarray*}
\mathcal{T}_{0}(d,u)&=&  \big(\alpha+\beta d^{2}+ \gamma
d^{4}\big)\Big[-2\varepsilon u+2U(d,u)\Big]\dfrac{1}{\tau_2(u)},\\
\mathcal{T}_{1}(d,u)&=&\big(\alpha+\beta d^{2}+ \gamma d^{4}\big)^{2}\dfrac{1}{\tau_2(u)^2}\mathcal{H}(d,u)^{2}
+\mathcal{R}(d,u)^{2}\\
&+&2\big(\alpha+\beta d^{2}+ \gamma d^{4}\big)\dfrac{1}{\tau_2(u)}\mathcal{H}(d,u)\mathcal{R}(d,u)\\
\mathcal{T}_{2}(d,u)&=& 2U(d,u)\mathcal{R}(d,u)+\frac{1}{2}\mathcal{R}(d,T^+)-\frac{1}{2}\mathcal{R}(d,T^-),\\
\mathcal{T}_{3}(d,u)&=&U(d,u)^{2}+ \mathcal{S}_2(U)-\varepsilon^{2}+\frac{d}{2}\mathcal{D}(d,u)
\end{eqnarray*}
Observing that $\mathcal{T}_{0}(d,u)$ is  analytic function beginnings with $d^3$. Using (1) of  lemma \ref{l6.6},  we obtain  that $d\mathcal{T}_{0}(d,u)\mathcal{H}(d,u)$ has property \textmd{G}.
 
Applying   lemma  \ref{l6.6}, we obtain
\begin{eqnarray*}
\mathcal{T}_{1}(d,u)=d^3\mathcal{D}_3(d,u)+\mathcal{O}\Big(d^3 e^{-\frac{\pi^2}{2d}}\Big)
\end{eqnarray*}
where $\mathcal{D}_{3}(d,u)$ has property \textmd{G}.

Using (5) of lemma \ref{l6.6} and  theorem \ref{t6.4}, we find  
\begin{eqnarray*}
\mathcal{T}_{2}(d,u)=d^3\mathcal{D}_4(d,u)+\mathcal{O}\Big(d^3 e^{-\frac{\pi^2}{2d}}\Big)
\end{eqnarray*}
where $\mathcal{D}_{4}(d,u)$ has property \textmd{G}. 

On the other hand, $\dfrac{1}{d}\mathcal{T}_{3}(d,u)$  has property \textmd{G}. 
This implies 
\begin{equation*}
\mathcal{T}(d,u)=\frac{1}{d}\mathcal{D}_{5}(d,u)+\mathcal{K}(d,u),
\end{equation*}
where $\mathcal{D}_{5}(d,u)$  has property \textmd{G} and 
\begin{equation}
\big|\mathcal{K}(d,u)\big|\leq Kd^{3}e^{\frac{-\pi^{2}}{2d}}\quad
\textrm{for}\quad (0<d<d_{0}, \,u_0<u<1).
\label{7.6} \end{equation}
 Altogether we have proved that 
\begin{equation*}
\mathcal{T}(d,u)=\frac{1}{d}\int_{0}^{\infty}e^{-\frac{t}{d}}\delta(t,u)dt+\mathcal{K}(d,u),
\end{equation*}
where $\delta(.,u)$ is analytic in $|t|<\frac{\pi^{2}}{2}$, it is
continuous on $[0,\frac{\pi^{2}}{2}[$ and $]\frac{\pi^{2}}{2},\infty[$ 
and has a limit as $t\to \frac{\pi^{2}}{2}$ for every $u_0<u<1$ and
satisfaies 
\begin{equation}
|\delta(t,u)|\leq K\,e^{Kt},
\label{7.7}\end{equation}
with  a constant  $K$ independent of $u$. Now if
$\delta(t,u)=\sum_{n=0}^{\infty}\delta_{n}(t)t^{n}$ is the power series
of $\delta(t,u)$ near $t=0$, Watson's lemma implies with (\ref{7.6}) that 
\begin{equation*}
d\mathcal{T}(d,u)\sim
\sum_{n=0}^{\infty}n\,!\delta_{n}d^{n+1}, \quad \textrm {as} \ d\searrow 0\  \textrm
{for every}\,\   u_0<u<1.
\end{equation*}
On the other hand, because of its definition 
\begin{eqnarray*}
\mathcal{T}(d,u)&\sim&\frac{1}{2}B(d,T^{+})-\frac{1}{2}B(d,T^{-})-\varepsilon^{2}+B(d,u)^{2}\\
&\sim & 0+0d+\cdot\cdot\cdot,
\end{eqnarray*}
because the formal series $B$ satisfy (\ref{5.1}). But this means that
all  $\delta_{n}$ must vanish and
$\delta(t,u)\equiv 0$ if \\ $|t|<\frac{\pi^{2}}{2}, \, u_0<u<1$. With (\ref{7.7}), thus we obtain for $u_0<u<1$ 
\begin{equation*}
\int_{0}^{\infty}e^{-\frac{t}{d}}|\delta(t,u)|dt\leq \int_{\frac{\pi^2}{2}}^{\infty}e^{-\frac{t}{d}}Ke^{Kt}dt\leq
Kde^{-\frac{\pi^{2}}{2d}},    \quad (0<d<d_{0}).
\end{equation*}
We obtain 
\begin{equation*}
\big|\mathcal{T}(d,u)\big|\leq Ke^{-\frac{\pi^{2}}{2d}}    \quad (0<d<d_{0}).
\end{equation*}
Finally, we have shown that there is a constant $K$ such that
\begin{equation*}
\Big|\mathcal{A}(d,T^{+})-\mathcal{A}(d,T^{-})-f\big(\varepsilon,\mathcal{A}(d,u)\big)\Big|=\frac{1}{\varepsilon}\big|\mathcal{T}(d,u)\big|\leq 
\frac{K}{d}e^{-\frac{\pi^{2}}{2d}},
\end{equation*}
for $0<d<d_{0}, \,0\leq u<1$.

\section{Distance between points of manifolds }
Let $y(t), t\in ]-1, \infty[$ an exact solution to the difference
equation (\ref{2.1}), so that the stable manifold $W^{+}_{s}$ is parametrized by
$t\mapsto \big(y(t),y(t+\varepsilon)\big)$. We know that
$y^{-}(t)=-y^{-}(-t)$ is also an exact solution of (\ref{2.1}) with
$y^{-}(t)\to 1 $ as $t\to \ -\infty $  .
Morever the instable manifold $W^{-}_{i}$ is parametrized by $t\mapsto
\big(y^{-}(t),y^{-}(t+\varepsilon)\big)$. We introduced the
quasi-solution $\mathcal{A}(d,u)$ in  previous section, then  
$\mathcal{A}^{-}(d,u)=-\mathcal{A}(d,-u)$ is a quasi-solution close to 
the exact solution $y^{-}(t)$. We denote by $\tilde{W}_{s}$,
$\tilde{W}_{i}$ the manifolds close to $W^{+}_{s}$  respectively  $W^{-}_{i}$  parametrized by  $t\mapsto
\big(\xi_{+}(t),\xi_{+}(t+\varepsilon)\big)$ respectively  $t\mapsto
\big(\xi_{-}(t),\xi_{-}(t+\varepsilon)\big)$, where
$\xi_{+}(t)=\mathcal{A}\big(d,u(t)\big)$ and $\xi_{-}(t)=\mathcal{A}^{-}\big(d,u(t)\big)$.   

In this section we will show that the
vertical distance between a point of the stable manifold and the
manifold $\tilde{W}_{s}$ is exponentially small as well as that
between $\tilde{W}_{i}$ and $W^{-}_{i}$; both are smaller, as we will show than the 
distance between  $\tilde{W}_{s}$ and  $\tilde{W}_{i}$. This eventually proves
that $W^{+}_{s}$ and  $W^{-}_{i}$ do not coincide.

In order to estimate the distance between some point $(x_{0},y_{0})$
on the stable manifold  $W^{+}_{s}$, we consider the sequence  $P_{n}=(x_{n},y_{n})$, defined by
$P_{n+1}=\phi(P_{n}),$ where 
\begin{eqnarray}
\phi\binom{x}{y}=\binom{y}{x+2\varepsilon(1-y^{2})}.
\label{8.1}\end{eqnarray}
 $P_{n}=(x_{n},y_{n})$ is  a  point of the stable manifold. We obtain
\begin{eqnarray}
P_{n+1}=\phi\binom{x_{n}}{y_{n}}=\binom{x_{n+1}}{y_{n+1}}=\binom{y_{n}}{x_{n}+2\varepsilon(1-y_{n}^{2})}
\label{8.2}\end{eqnarray}
 Since $\Phi=\phi\circ\phi$, where $\Phi$ is defined in (\ref{1.3}), and  the stable manifold   $W_{s}^{+}$ is invariant by the
 application $\Phi$, then the stable  manifold   $W_{s}^{+}$  remains invariant by  $\phi$. Then theorem \ref{t7.1} implies 
\begin{eqnarray}
\xi_{+}(t+\varepsilon)=\xi_{+}(t-\varepsilon)+2\varepsilon\big(1-\xi_{+}(t)^{2}\big)+\frac{1}{\varepsilon}\eta(t)
e^{-\frac{\pi^{2}}{2\varepsilon}}
\label{8.3}\end{eqnarray}
where
$\eta(t)$ is bounded.
There is a sequence $t_{n}\in ]-1,+\infty[$ such that 
\begin{eqnarray}
\xi_{+}(t_{n}-\varepsilon)=x_{n}=\mathcal{A}\Big(d,T^{-}\big(d,u(t_{n})\big)\Big)=\mathcal{A}\Big(d,T^{-}\big(d,U_{\varepsilon}(n)\big)\Big)
\label{8.4}\end{eqnarray}
\begin{eqnarray}
\xi_{+}(t_{n})&=&\mathcal{A}\big(d
,U_{\varepsilon}(n)\big)\nonumber\\ 
\xi_{+}(t_{n}+\varepsilon)&=&\mathcal{A}\Big(d,T^{+}\big(d,U_{\varepsilon}(n)\big)\Big)
\label{8.5}\end{eqnarray}
where $U_{\varepsilon}(n)=u(t_{n})$. 
 The relation (\ref{8.4}) implies that the point $(x_{n}, b_{n})$, where $b_{n}=\xi_{+}(t_{n})$,  the vertical projection
 of the point $(x_{n}, y_{n})$ on the curve $t\mapsto
 \big(\xi_{+}(t),\xi_{+}(t+\varepsilon)\big)$. But, according to
 (\ref{8.2}), (\ref{8.4}), we have 
 \begin{eqnarray}
\xi_{+}(t_{n+1}-\varepsilon)=x_{n+1}=y_{n}.
\label{8.6}\end{eqnarray}
\begin{figure}[!t]
  \centering
  
   
  \includegraphics[width=.4\linewidth]{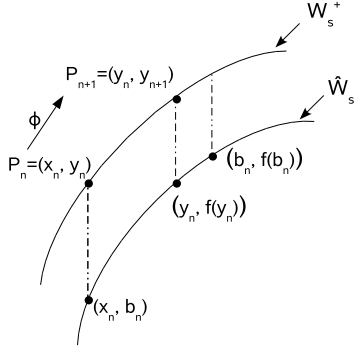}
  
  \caption{\small{The vertical distance between a point of the stable manifold and the
approximation defined using the quasi-solution $\mathcal{A}$
      for the logistic equation .}}
\end{figure}
Now consider the error quantities 
\begin{eqnarray}
d_{\varepsilon}(n)=y_{n}-\mathcal{A}\big(d,U_{\varepsilon}(n)\big)=y_{n}-b_{n},
\label{8.7}\end{eqnarray}
$d_{\varepsilon}(n)$ denotes the vertical distance between a  point $P_{n}$ and the 
manifold $\Hat{W}_{s}$. The manifold $\tilde{W}^{+}_{s}$ is  parameterized by the curve\\
$t\mapsto \big(\xi_{+}(t),\xi_{+}(t+\varepsilon)\big)$.
 We obtain
\begin{eqnarray}
d_{\varepsilon}(n+1)=y_{n+1}-f(x_{n+1})=y_{n+1}-f(y_{n})=y_{n+1}-f\big(d_{\varepsilon}(n)+b_{n}\big),
\label{8.8}\end{eqnarray}
where $f(x)=\xi_{+}\Big(\xi_{+}^{-1}(x)+\varepsilon\Big)$. The point
$(y_{n},f\big(y_{n})\big)$ is the  vertical projection
 of the point $(x_{n+1}, y_{n+1})$ on the curve $t\mapsto
 \big(\xi_{+}(t),\xi_{+}(t+\varepsilon)\big)$\quad (See FIG.4).

If we use Taylor expainson, (\ref{8.8}) implies
\begin{eqnarray}
d_{\varepsilon}(n+1)=y_{n+1}-f(b_{n})-f'(\theta_{n})d_{\varepsilon}(n),\qquad 
b_{n}<\theta_{n}<b_{n}+d_{\varepsilon}(n).
\label{8.9}\end{eqnarray}
Then according to the definitions of $f$ and $b_{n}$, we have
 $f(b_{n})=\xi_{+}(t_{n}+\varepsilon)$, with (\ref{8.3}) and (\ref{8.4}), we obtain 
\begin{eqnarray}
f(b_{n})=x_{n}+2\varepsilon(1-b_{n}^{2})+\frac{1}{\varepsilon}\eta(t_{n})e^{-\frac{\pi^{2}}{2\varepsilon}}.
\label{8.10}\end{eqnarray}
Together with (\ref{8.9}) and (\ref{8.2}) this implies

\begin{eqnarray}
d_{\varepsilon}(n+1)=-f'(\theta_{n})d_{\varepsilon}(n)+2\varepsilon(b^{2}_{n}-y^{2}_{n})-\frac{1}{\varepsilon}\eta(t_{n})e^{-\frac{\pi^{2}}{2\varepsilon}},
\label{8.11}\end{eqnarray}
Using (\ref{8.7}) we obtain 
\begin{eqnarray}
d_{\varepsilon}(n+1)=-\Big(f'(\theta_{n})+2\varepsilon(y_{n}+b_{n})\Big)d_{\varepsilon}(n)-\frac{1}{\varepsilon}\eta(t_{n})e^{-\frac{\pi^{2}}{2\varepsilon}}.
\label{8.12}\end{eqnarray}
Now 
\begin{equation*}
f'(\theta_{n})=\frac{\xi'_{+}(z_{n}+\varepsilon)}{\xi'_{+}(z_{n})},\quad
\text{where}\quad z_{n}=\xi_{+}^{-1}(\theta_{n}),
\end{equation*}
then 
\begin{equation*}
f'(\theta_{n})=1-2u(z_{n})d+\big(3u(z_{n})^{2}-1\big)d^{2}+O(d^{3}).
\end{equation*}
On the other hand $b_{n}<\theta_{n}<y_{n}$  implies 
\begin{equation*}
0<b_{n}<\xi(z_{n})<y_{n}
\end{equation*}
and thus 
\begin{eqnarray}
f'(\theta_{n})+2\varepsilon(y_{n}+b_{n})&>&
1+2\varepsilon\Big(\xi_{+}(z_{n})+b_{n}\Big)-2u(z_{n})d\nonumber \\
&+&\Big(3u(z_{n})^{2}-1\Big)d^{2}+O(d^{3}).
\label{8.13}\end{eqnarray}
We know that  
\begin{equation*}
\xi_{+}(z_{n})=u(z_{n})+\Big(u(z_{n})-u(z_{n})^{3}\Big)d^{2}+O(d^{4}),
\end{equation*}
this with (\ref{2.8}) implies 
\begin{eqnarray*}
f'(\theta_{n})+2\varepsilon(y_{n}+b_{n})>
1+2\varepsilon b_{n}+O(d^{2}).
\end{eqnarray*}
As $ b_{n}>0$,
then for sufficiently small $d$,
\begin{eqnarray*}
f'(\theta_{n})+2\varepsilon(y_{n}+b_{n})>1.
\end{eqnarray*}
This implies
\begin{eqnarray*}
f'(\theta_{n})+2\varepsilon(y_{n}+b_{n})=
1+g_{\varepsilon}\big(u(z_{n})\big), \quad \text{with some}\quad
g_{\varepsilon}\big(u(z_{n})\big)\geq \varepsilon c,
\end{eqnarray*}
 $g_{\varepsilon}\big(u(z_{n})\big)=O(\varepsilon)$ and $c$ is 
 a positive constant.  
With (\ref{8.12}), this implies 
\begin{eqnarray}
d_{\varepsilon}(n+1)=-\Big(1+g_{\varepsilon}\big(u(z_{n})\big)\Big)d_{\varepsilon}(n)-\frac{1}{\varepsilon}\eta(t_{n})e^{-\frac{\pi^{2}}{2\varepsilon}}.
\label{8.14}\end{eqnarray}

As  $d_{\varepsilon}(n)\to 0$ as  $n\to +\infty$, this implies that
\begin{equation}
|d_{\varepsilon}(n)|\leq \frac{1}{\varepsilon}e^{-\frac{\pi^{2}}{2\varepsilon}} 
\sum^{\infty}_{k=
  n}\bigg |\frac{\eta(s_{k})}{\prod_{i}^{i=k}\Big(1+g_{\varepsilon}\big(u(z_{k})\big)\Big)}\bigg |,
\label{8.15}\end{equation}
provided this series converges. This follows from
\begin{eqnarray*}
\sum^{\infty}_{k=n}\bigg |\frac{\eta(s_{k})}{\prod_{i}^{i=k}\Big(1+g_{\varepsilon}\big(u(z_{k})\big)\Big)}\bigg |
\leq K_{1}\sum^{+\infty}_{k=0}(1+\varepsilon c\big)^{-k},
\end{eqnarray*}
where we used that $\eta$ is bounded. 
This imply that there is a constant $K>0$ such that 
\begin{equation*}
\sum^{\infty}_{k=n}\bigg |\frac{\eta(s_{k})}{\prod_{i}^{i=k}\Big(1+g_{\varepsilon}\big(u(z_{k})\big)\Big)}\bigg |
\leq  \frac{K}{\varepsilon}.
\end{equation*}
 With (\ref{8.5}) implies
\begin{eqnarray*}
d_{\varepsilon}(n)=O\bigg( \frac{1}{\varepsilon^{2}}\exp\Big(-\frac{\pi^{2}}{2\varepsilon}\Big)\bigg),
\end{eqnarray*}
and in particular that the vertical distance $d_{\varepsilon}(0)$ of
$(x_{0},y_{0})$ to $\tilde{W}_{s}$ is 
\begin{eqnarray}
d_{\varepsilon}(0)=O\bigg( \frac{1}{\varepsilon^{2}}\exp\Big(-\frac{\pi^{2}}{2\varepsilon}\Big)\bigg).
\label{8.16}\end{eqnarray}

 Now we will calculate the vertical distance 
between the twomanifolds  $\tilde{W}_{s}$ and $\tilde{W}_{i}$. To do this, we will need to extend both
quasi-solutions in an adequate way.
We define 
\begin{eqnarray}
\tilde{\xi}_{+}(t)&=&\mathcal{A}\big(d,u(t)\big), \quad \textrm{for} \quad
-\frac{1}{2}\leq u\leq 1, \nonumber\\
\tilde{\xi}_{-}(t)&=&\mathcal{A}^{-}\big(d,u(t)\big), \ \, \textrm{for} \quad
-1\leq u\leq \frac{1}{2}.
\label{8.17}\end{eqnarray}
First, we  evaluate the  quantities
\begin{equation}
 d_{\varepsilon}\Big(\tilde{\xi}_{+}(t),\tilde{\xi}_{-}(t)\Big)=\mathcal{A}(d,u)-\mathcal{A}^{-}(d,u),  \ \textrm{for} \
-\frac{1}{2}\leq u\leq \frac{1}{2}
\label{8.18}\end{equation}
using (\ref{7.2}) and (\ref{7.4}) 
\begin{eqnarray}
 d_{\varepsilon}\Big(\tilde{\xi}_{+}(t),\tilde{\xi}_{-}(t)\Big)&=&\tilde{\xi}_{+}(t)-\tilde{\xi}_{-}(t),\\ 
&=&\frac{1}{\varepsilon}(\alpha+\beta d^{2}+\gamma d^{4})\big(\mathcal{\tilde{H}}(d,u)-\tilde{\mathcal{H}}^{-}(d,u)\big),\nonumber
\label{8.19}\end{eqnarray}
where 
\begin{eqnarray}
\tilde{\mathcal{H}}(d,u)&=&\int^{\infty}_{0}\, e^{-\frac{s}{d}}h(s,u)\,ds, \quad \textrm{for} \quad
-\frac{1}{2}\leq u\leq \frac{1}{2}, \nonumber\\
\tilde{\mathcal{H}}^{-}(d,u)&=&-\int^{\infty}_{0}\, e^{-\frac{s}{d}}h(s,-u)\,ds, \quad \textrm{for} \quad
-\frac{1}{2}\leq u\leq \frac{1}{2}.
\label{8.20}\end{eqnarray}
Modifying the path of integral, we obtain 
\begin{eqnarray*}
\tilde{\mathcal{H}}(d,u)=\int_{\Gamma}\, e^{-\frac{s}{d}}h(s,u)\,ds&+&
\sum_{Im(s_{k}(u))<0}2\pi i
Res\Big(e^{-\frac{s}{d}}h(s,u),s_{k}(t)\Big)\\
&-&\sum_{Im(s_{k}(u))>0}2\pi i
Res\Big(e^{-\frac{s}{d}}h(s,u),s_{k}(t)\Big),
\end{eqnarray*}
where $s_{k}(t)= \frac{\pi^{2}}{2}\pm
\frac{d\pi\,t}{\varepsilon}\,i+k\pi^{2}$ \,for $k\ge 0$,  $\Gamma$ is
the sum of two paths $\Gamma_{1}$ and  $\Gamma_{2}$. Here
$\Gamma_{1}$ consist of two segments of the higher half-plane, one of
these two segments is parallel to the axis $y=0$ and begins at the  point $(\frac{\pi^{2}}{4},1)$, the other ends
in the point (0,0), and  the path $\Gamma_{2}$  is the symmetry of
$\Gamma_{1}$ relatively to the axis $y=0$.
 
This implies
\begin{eqnarray*}
\tilde{\mathcal{H}}(d,u)=\int_{\Gamma}\, e^{-\frac{s}{d}}h(s,u)\,ds&+&
 \frac{\pi e^{\frac{\pi\,t}{\varepsilon}i}}{2d^{2}(1-u^{2})}\sum_{k=0}^{\infty}
\exp\bigg(-\frac{(k+1)\pi^{2}}{2d}\bigg)\\
&+&\frac{\pi e^{-\frac{\pi\,t}{\varepsilon}i}}{2d^{2}(1-u^{2})}\sum_{k=0}^{\infty}
\exp\bigg(-\frac{(k+1)\pi^{2}}{2d}\bigg)
\end{eqnarray*}
Therefore
\begin{eqnarray}
\tilde{\mathcal{H}}(d,u)=\int_{\Gamma}\, e^{-\frac{s}{d}}h(s,u)\,ds+
 \frac{\pi\cos(\frac{\pi\,t}{\varepsilon})}{d^{2}(1-u^{2})}\sum_{k=0}^{\infty}
\exp\bigg(-\frac{(k+1)\pi^{2}}{2d}\bigg).
\label{8.21}\end{eqnarray}
Similary for\ $\tilde{\mathcal{H}}^{-}(d,u)$
\begin{eqnarray*}
\tilde{\mathcal{H}}^{-}(d,u)=\int_{\Gamma}\, e^{-\frac{s}{d}}h(s,u)\,ds-
 \frac{\pi\cos(\frac{\pi\,t}{\varepsilon})}{d^{2}(1-u^{2})}\sum_{k=0}^{\infty}
\exp\bigg(-\frac{(k+1)\pi^{2}}{2d}\bigg).
\end{eqnarray*}
 Therefore 
\begin{eqnarray*}
  d_{\varepsilon}\Big(\tilde{\xi}_{+}(t),\tilde{\xi}_{-}(t)\Big)=
 \frac{2\pi(\alpha+\beta d^{2}+\gamma d^{4})\cos(\frac{\pi\,t}{\varepsilon})}{\varepsilon^{3} \,(1-u^{2})}
\exp\bigg(-\frac{\pi^{2}}{2d}\bigg)+\mathcal{O}\bigg(\frac{\cos(\frac{\pi\,t}{\varepsilon})}{\varepsilon^{3}}e^{-\frac{3\pi^{2}}{2d}}\bigg), 
\end{eqnarray*}
for  $-\text{arctanh}(\frac{1}{2})<t<\text{arctanh}(\frac{1}{2})$.
Consequently
\begin{eqnarray}
  d_{\varepsilon}\Big(\tilde{\xi}_{+}(t),\tilde{\xi}_{-}(t)\Big)=
  \frac{2\pi\alpha\cos(\frac{\pi\,t}{\varepsilon})}{\varepsilon^{3}\big(1-\tanh(\frac{d}{\varepsilon}t)^{2}\big)}e^{-\frac{\pi^{2}}{2\varepsilon}}+\mathcal{O}\bigg(\frac{1}{\varepsilon}\cos(\frac{\pi\,t}{\varepsilon})e^{-\frac{\pi^{2}}{2d}}\bigg)
\label{8.22}\end{eqnarray}

 Now we want to estimate the  vertical distance  between a point on the
 quasimanifold $\tilde{W}_{s}$  and its vertical projection on the  quasimanifold   $\tilde{W}_{i}$.
 For this purpose  let $t$ such that $-\frac{1}{2}\leq u(t)\leq
 \frac{1}{2}$. Then the  point
 $\Big(\tilde{\xi}_{+}(t),\tilde{\xi}_{+}(t+\varepsilon)\Big)$ is on
 the  quasimanifold $\tilde{W}_{s}$. We suppose that
 $\Big(\tilde{\xi}_{+}(t),\tilde{\xi}^{-}(t_{1}+\varepsilon)\Big)$,
 where  $-\frac{1}{2}\leq u(t_{1})\leq \frac{1}{2}$,  is its vertical
 projection on the  quasimanifold   $\tilde{W}_{i}$, then the vertical distance  between these two points is 
\begin{eqnarray}
dist_{v}(t)=\tilde{\xi}_{+}(t+\varepsilon)-\tilde{\xi}_{-}(t_{1}+\varepsilon)=\tilde{\xi}_{+}(t+\varepsilon)-\tilde{f}\big(\tilde{\xi}_{+}(t)\big),
\label{8.23}\end{eqnarray}
where
$\tilde{f}(x)=\tilde{\xi}_{-}\Big(\tilde{\xi}_{-}^{-1}(x)+\varepsilon\Big)$, 
$\tilde{\xi}_{-}^{-1}$ the inverse function of $\tilde{\xi}_{-}$.
Using (\ref{8.22}) and (\ref{8.19}) we obtain 
\begin{eqnarray*}
dist_{v}(t)=\tilde{\xi}_{+}(t+\varepsilon)-\tilde{f}\big(\tilde{\xi}_{-}(t)+e_{\varepsilon}(t)\big),
\end{eqnarray*}
where\ \ $e_{\varepsilon}(t)= d_{\varepsilon}\Big(\tilde{\xi}_{+}(t),\tilde{\xi}_{-}(t)\Big)$.
We use Taylor expansion and obtain 
\begin{eqnarray}
dist_{v}(t)&=&\tilde{\xi}_{+}(t+\varepsilon)-\tilde{f}\big(\tilde{\xi}_{-}(t)\big)-e_{\varepsilon}(t)\cdot 
\tilde{f}'(\theta),\nonumber\\
&=&\tilde{\xi}_{+}(t+\varepsilon)-\tilde{\xi}_{-}(t+\varepsilon)-e_{\varepsilon}(t)\cdot 
\tilde{f}'(\theta),
\label{8.24}\end{eqnarray}
where \  \
$\tilde{\xi}_{-}(t)<\theta<\tilde{\xi}_{-}(t)+e_{\varepsilon}(t)$.
With (\ref{8.19}) this implies 
\begin{eqnarray}
dist_{v}(t)&=&e_{\varepsilon}(t+\varepsilon)-e_{\varepsilon}(t)\cdot\tilde{f}'(\theta).
\label{8.25}\end{eqnarray}
Since\ 
$e_{\varepsilon}(t+\varepsilon)=-e_{\varepsilon}(t)+\mathcal{O}\Big(\varepsilon\cdot 
e_{\varepsilon}(t)\Big)$ 
and
$\tilde{f}'(\theta)=\dfrac{\tilde{\xi}'_{-}\Big(\tilde{\xi}^{-1}_{-}(\theta)+\varepsilon\Big)}{\tilde{\xi}'_{-}\Big(\tilde{\xi}^{-1}_{-}(\theta)\Big)}=1+\mathcal{O}(\varepsilon)$,
we obtain 
\begin{eqnarray}
dist_{v}(t)&=&\big(-2+\mathcal{O}(\varepsilon)\big)e_{\varepsilon}(t).
\label{8.26}\end{eqnarray}
Consequently
\begin{equation}
dist_{v}(t)\sim- \frac{4\pi\alpha\cos(\frac{\pi\,t}{\varepsilon})}{\varepsilon^{3}\big(1-\tanh(\frac{d}{\varepsilon}t)^{2}\big)}e^{-\frac{\pi^{2}}{2\varepsilon}}.
\label{8.27}\end{equation}
Combining (\ref{8.16}) and (\ref{8.27}), we conclude
\begin{equation*} 
Dist_{\varepsilon}\big(w^{+}_{s}(t),W^{-}_{i}\big)=\frac{4\pi\alpha\cos(\frac{\pi}{\varepsilon}t+\pi)}{\varepsilon^{3}\big(1-\tanh(\frac{d}{\varepsilon}t)\big)}e^{\frac{-\pi^{2}}{2\varepsilon}}+O\Bigg(\frac{1}{\varepsilon^{2}}e^{\frac{-\pi^{2}}{2\varepsilon}}\Bigg),\quad \text{as} \,\,\,\varepsilon\searrow 0,
\end{equation*}

\section{An estimate for the main asymptotic coefficients}
In this section we show that the number $\alpha$ in the main theorem of the
introduction is not zero  and we want to find bounds for
$\alpha_{n}$ and will then obtain estimats for $\alpha$  using the
first $\alpha_{n}$ which can be computed explicitly. We recall that 
\begin{equation}
\alpha=\frac{1}{\pi}\sum_{n=7}^{\infty}\alpha_{n}
\label{9.1}\end{equation}
where $\alpha_{n}(n-1)\,!\big(\frac{i}{\pi}\big)^{n-1}u^{n}$ is the
leading  term of the coefficient $ \{\mathcal{C}(G)\}_{n}(u)$ of
$d^{n}$ in  $ \mathcal{C}(G)(d,u)$, where $G=Q(d,u).F(d,u),\,
F(d,u)=B(d,u)-U(d,u)$ \,(cf.(\ref{5.2}) \,and (\ref{5.5})).

For a polynomial series $X(d,u)=\sum_{n=0}^{\infty}X_{n}(u)d^{n}$
where the
degrees of $X_{n}(u)$ do not exceed $n$, we can write
$X_{n}(u)=\sum_{k=0}^{n}x_{nk}u^{n}$.
Then we denote by $\widehat{X}(z)=\sum_{n=0}^{\infty}x_{nn}z^{n}$.

The mapping $X(d,u)\mapsto \widehat{X}(z)$ extracts the leading terms
of the series $X(d,u)$. It is compatible with addition 
and multiplication.
If we define the following operator  $\widehat{D}=-z^{2}\frac{d}{dz}$, we also find
\begin{eqnarray*}
\widehat{D}\widehat{X}(z)&=&\widehat{DX}(z)=-z^{2}\widehat{X}'(z)\\
 \widehat{\mathcal{C}(X)}(z)&=&\cosh\bigg( \frac{\widehat{D}}{2}\bigg)
 \widehat{X}(z):=\frac{1}{2}\Bigg(\widehat{X}\Big(\frac{2z}{2+z}\Big)+\widehat{X}\Big(\frac{2z}{2-z}\Big)\Bigg)\\
\widehat{\mathcal{S}(X)}(z)&=&\sinh\bigg( \frac{\widehat{D}}{2}\bigg)
 \widehat{X}(z):=\frac{1}{2}\Bigg(\widehat{X}\Big(\frac{2z}{2+z}\Big)-\widehat{X}\Big(\frac{2z}{2-z}\Big)\Bigg)\\
\widehat{\mathcal{C}_{2}(X)}(z)&=&\cosh( \widehat{D})
 \widehat{X}(z):=\frac{1}{2}\Bigg(\widehat{X}\Big(\frac{z}{1+z}\Big)+\widehat{X}\Big(\frac{z}{1-z}\Big)\Bigg)\\
\widehat{\mathcal{S}_{2}(X)}(z)&=&\sinh( \widehat{D})
 \widehat{X}(z):=\frac{1}{2}\Bigg(\widehat{X}\Big(\frac{z}{1+z}\Big)-\widehat{X}\Big(\frac{z}{1-z}\Big)\Bigg)
\end{eqnarray*}
where the operators $\mathcal{S},\mathcal{C},\mathcal{S}_{2},\mathcal{C}_{2}$
are defined in (\ref{2.18}). 
If we apply the mapping $X(d,u)\mapsto \widehat{X}(z)$ to the
equation (\ref{5.13}) , we obtain
\begin{equation}
\widehat{\mathcal{L}}\bigg(\frac{\widehat{E}}{\widehat{Q}_{1}}\bigg)(z)=\widehat{\mathcal{R}}\big(\widehat{E}\big)(z),
\label{9.2}\end{equation}
where
\begin{eqnarray}
\widehat{\mathcal{L}}(x)&:=&-\frac{\widehat{V}_{1}}{z^{2}}\bigg(\widehat{\mathcal{S}}(x)-\frac{\widehat{\mathcal{S}}(z^{2})}{\widehat{\mathcal{C}}(z^{2})}\widehat{\mathcal{C}}(x)\bigg)\nonumber\\
\widehat{\mathcal{R}}(x)&:=&
-\frac{\widehat{\mu}}{z^{2}}+\frac{\widehat{W}_{5}}{z^{2}}\widehat{\mathcal{C}}\bigg(\frac{x}{\widehat{Q}_{1}}\bigg)-\frac{\widehat{W}_{1}}{2z^{2}\widehat{Q}}\mathcal{Y}+\frac{\widehat{W}_{4}}{2z^{2}\widehat{Q}^{2}}\widehat{\mathcal{C}}\bigg(\frac{x}{\widehat{Q}_{1}}\bigg)\mathcal{Y}^{2}\nonumber\\
&-& \frac{\widehat{W}_{2}}{2z^{2}}\widehat{\mathcal{C}}_{2}\bigg(\frac{\mathcal{Y}}{\widehat{Q}}\bigg)\nonumber,\\
\widehat{E}&=&\widehat{\mathcal{C}(G)}=\sum_{n=0\, \text{odd}}^{\infty}\alpha_{n}(n-1)\,!\bigg(\frac{i}{\pi}\bigg)^{n-1}z^{n},\nonumber\\
\mathcal{Y}&=&\widehat{\mathcal{C}}^{-1}(x),\nonumber\\
\widehat{Q}&=&-z^{2}+z^{4}-\frac{13}{6}z^{6}+\frac{47}{18}z^{8},\nonumber\\
\widehat{Q}_{1}&=&-z^{2}-\frac{3}{2}z^{4},\nonumber\\
\widehat{V}_{1}&=& 4z^{2}\bigg[-192-240z^{2}+508z^{4}+1995z^{6}+3418z^{8}\nonumber\\
&+&4389z^{10}+217z^{12}\bigg]\bigg/\bigg[3(z^{2}-4)^{4}\bigg],\nonumber\\
\widehat{W}_{2}&=&z^{11}\bigg[591-15457z^{2}+21610z^{4}+59194z^{6}+36583z^{8}\nonumber\\
&+&2759z^{10}\bigg]\bigg/\bigg[27(z^{2}-1)^{8}\bigg],\nonumber\\
\widehat{W}_{1}&=& z^{11}\bigg[2385-32832z^{2}-11790z^{4}+231917z^{6}-116267z^{8}+720682z^{10}\nonumber\\
&+&1018936z^{12}+290089z^{14}+2480z^{16}\bigg]\bigg/\bigg[81(z^{2}-1)^{8}\bigg],\nonumber\\
\widehat{W}_{4}&=&\widehat{V}\widehat{\mathcal{C}}_{2}(\widehat{Q})+
\widehat{W}\widehat{\mathcal{S}}_{2}(\widehat{Q})\nonumber\\
&=&-z^{2}\bigg[54-432z^{2}+1449z^{4}-2724z^{6}+4449z^{8}\nonumber 
-14719z^{10}-12857z^{12}\nonumber\\
&+&95146z^{14}+110905z^{16}+29041z^{18}+248z^{20}\bigg]\bigg/\bigg[54(z^{2}-1)^{8}\bigg],\nonumber\\
\widehat{W}_{5}&=&\widehat{W}_{3}+\widehat{P}\widehat{V}_{3},\nonumber\\
&=&-8z^{5}\bigg[-1152+1168z^{2}+9744z^{4}+18901z^{6}+14664z^{8}\nonumber\\
&+&1240z^{10}\bigg]\bigg/\bigg[3(z^{2}-4)^{4}(4+z^{2})\bigg],\nonumber
\end{eqnarray}
\begin{eqnarray}
\widehat{\mu}&=& z^{10}\bigg[-527000z^{30}-1988300z^{28}+14763570z^{26}-25488145z^{24}\nonumber\\
&+&11808105z^{22}+11381027z^{20}-17065562z^{18}-15479497z^{16}-15777878z^{14}\nonumber\\ 
&-&21499321z^{12}+14844222z^{10}+4847178z^{8}-1419003z^{6}-541026z^{4}\nonumber\\
&+&320490z^{2}-50220\bigg]\bigg/\bigg[972(z^{2}-1)^{12}\bigg],
\label{9.3}\end{eqnarray}
with $x=x(z)$ \,is a power series.

Now for integers $m> p\geq -1$,\, $m\geq 0$, we introduce the set
\begin{eqnarray*}
\mathcal{M}(m,p)=\bigg\{\summ{n}{m}{\infty} x_{n}z^{n}\backslash\,\,\, x_{n}\big/\big[
(n-p)\,!\pi^{-n}\big],n=m,m+2...\,\,\text{is bounded}\bigg\}
\end{eqnarray*}
$\mathcal{M}(m,p)$ \, is a Banach space with the norm
\begin{eqnarray}
|x|_{m,p}=\sup\big\{\, |x_{n}|\big/\big[
(n-p)\,!\pi^{-n}\big],n\geqq m, n\equiv m\mod 2\big\}.
\label{9.4}\end{eqnarray}

Here the symbol $\summ{n}{m}{\infty}$ for the summation means that
only summands with $ n\equiv m \mod 2$ are considered.
\begin{propo} Let m,p,q be integers and $m\geq p\geq -1,m\geq q\geq -1,m \geq 0 $.
  The norms (\ref{9.4}) have the following property:
\begin{enumerate}
\item If \, $q<p\le m$ then  \,  $\mathcal{M}(m,p)\subset\mathcal{M}(m,q)$  and 
 \begin{equation}
|x|_{m,q}\leq \frac{(m-p)\,!}{(m-q)\,!}|x|_{m,p}\quad \text{for}\,\,\, x\in\mathcal{M}(m,p)
 \label{9.5}\end{equation}
\item If $x\in\mathcal{M}(m,p)$ then $zx\in\mathcal{M}(m+1,p+1)$ and
\begin{equation}
|zx|_{m+1,p+1}=\pi|x|_{m,p}
 \label{9.6}\end{equation}

\item $\widehat{D}:\mathcal{M}(m,p) \mapsto\mathcal{M}(m+1,p)$ and 
\begin{equation}
\pi|x|_{m,p}\leq|\widehat{D}(x)|_{m+1,p}\leq\pi\frac{m}{m-p+1}|x|_{m,p}
\label{9.7} \end{equation}
\end{enumerate}
\label{p9.1}\end{propo}

Now we consider an even power series
$f(z)=\sum_{j=0}^{\infty}*f_{j}z^{j}$ having a radius of convergence
greater than $\pi$. Then $f(\widehat{D})$ maps $\mathcal{M}(m,p)$ into 
itself and 
\begin{equation}
|f(\widehat{D})x|_{m,p}\leq \|f\|_{m,p}|x|_{m,p},\quad \text{for} \,
x\in \mathcal{M}(m,p)
\label{9.8}\end{equation}
where 
\begin{eqnarray}
\|f\|_{m,p}&=&\sup\bigg\{\summ{j}{0}{n-m}|f_{j}|\pi^{j}\prod_{k=1}^{p-1}\frac{n-k}{n-j-k}\backslash 
n\geq m, n\equiv m \mod 2\bigg\}\,\text{if}\ p\geq 1,\nonumber\\
\|f\|_{m,p}&=&\sup\bigg\{\summ{j}{0}{\infty}|f_{j}|\pi^{j}\backslash 
n\geq m, n\equiv m \mod 2\bigg\}\,\text{if}\ p\leq 0.
\label{9.9}\end{eqnarray}

 We have also
\begin{equation}
|f.x|_{m,p}\leq M_{m,p}(f)|x|_{m,p},\quad \text{for} \,
x\in \mathcal{M}(m,p),
\label{9.10}\end{equation}
where
\begin{eqnarray}
M_{m,p}(f)=\sup\bigg\{\summ{j}{0}{n-m}|f_{j}|\pi^{j}\frac{(n-p-j)\,!}{(n-p)\,!}\backslash 
n\geq m, n\equiv m \mod 2\bigg\}\,.
\label{9.11}\end{eqnarray}
We conclude that the inverse of $\widehat{\mathcal{C}}$ maps $
\mathcal{M}(m,p)$ with $p\geq 3$ into   $\mathcal{M}(m,1)$ and 
\begin{eqnarray}
|\widehat{\mathcal{C}}^{-1}(x)|_{m,1}\leq
c_{m,p}|x|_{m,p},\ \text{with}\ c_{m,p}=\frac{4}{\pi}\summ{n}{m}{\infty}\frac{(n-p)\,!}{(n-1)\,!}.
\label{9.12}\end{eqnarray}

Let  $x\in \mathcal{M}(m,p)$ and $y\in \mathcal{M}(m,q)$ and assume
that $m+q\leq n+p$. Then we obtain $x.y \in\mathcal{M}(n+m,q+m)$  and 
\begin{eqnarray}
|x.y|_{n+m,q+m}\leq  \alpha_{m-p,n-q}|x|_{m,p}|y|_{n,q},
\label{9.13}\end{eqnarray}
with $$\alpha_{r,s}=\sup\Bigg\{\summ{j}{r}{N-s}\frac{j\,!(N-j)\,!}{(N-r)\,!}\backslash 
  N\geq r+s,N-r-s \,\text{even}\Bigg\}.$$
Finally we need estimates for the inverse of the linear operator
$\widehat{\mathcal{L}}$ in (\ref{9.2}).  $ \widehat{\mathcal{L}}$ maps $
\mathcal{M}(m-1,q)$ onto   $\mathcal{M}(m,q)$ and is one to one.Now 
let   $x\in \mathcal{M}(m-1,p)$  and $w= \widehat{\mathcal{L}}(x)$. If 
we proceed similar to the proof of theorem \ref{t4.6}, we find that 
\begin{eqnarray}
 \widehat{D}y&=&f(
 \widehat{D})\bigg[-\frac{(4+z^{2})}{4\widehat{V}_{1}}w\bigg],\nonumber\\
x&=& z^{2}y.
\label{9.14}\end{eqnarray}
where $f(z)=z/\sinh(z/2)$.\\

We can apply (\ref{9.6})..(\ref{9.13})  and obtain for $m\geq \max(p+1,4)$
\begin{eqnarray}
&&|\widehat{\mathcal{L}}^{-1}(w)|_{m-1,p}\leq L_{m,p}|w|_{m,p} \quad
\text{for}  \, w \in\mathcal{M}(m,p),
\label{9.15}\end{eqnarray}
where
\begin{eqnarray}
L_{m,p}=\frac{1}{\pi}\bigg\|\frac{z}{\sinh(z/2)}\bigg\|_{m-2,p-2}M_{m,p}\bigg(\frac{(4+z^{2})z^{2}}{4\widehat{V}_{1}}\bigg).
\label{9.16}\end{eqnarray}

Now we are in a position to use (\ref{9.2}) to estimate $\widehat{E}$. We
use that $\widehat{E}$ can be computed recursively from  (\ref{9.2}). If $x$ is an odd 
power series and $x\equiv\widehat{E}\mod z^{n} $ with $n\geq 9$ then
via $ \widehat{\mathcal{C}}^{-1}(x)\equiv\widehat{G}\mod z^{n}$ we
first obtain
$\widehat{\mathcal{R}}(x)\equiv\widehat{{\mathcal{R}}}(\widehat{E})\mod
z^{n+1}$ and then
$\widehat{\mathcal{L}}^{-1}\Big(\widehat{\mathcal{R}}(x)\Big)\equiv\frac{\widehat{E}}{\widehat{Q}_{1}}\mod
z^{n+1}$.Thus if $x_{0}\equiv 0 \mod z^{9}$ is arbitrary and we let
\begin{eqnarray}
x_{n+1}:=\widehat{Q}_{1}\widehat{\mathcal{L}}^{-1}\Big(\widehat{\mathcal{R}}(x_{n})\Big)
 \quad\text{for}\, n=0...
\label{9.17}\end{eqnarray} 
then $x_{n}\equiv \widehat{E}\mod z^{9+2n}$ for all
$n$. 
We begin with the polynomial $x_{0}$ of degree $\leq 79$ with
$x_{0}\equiv\widehat{E}\mod z^{81}$. We computed it explicitly and
found that 
\begin{eqnarray}
|x_{0}|_{9,7}\leq K_{0}:=308027.359777894414.
\label{9.18}\end{eqnarray}

Then we consider the sequence $x_{n}$  of polynomial series defined  
in (\ref{9.17}), which now satisfies  $x_{n}\equiv \widehat{E}\mod
 z^{81+2n}$. We will show 
\begin{eqnarray}
|x_{n}|_{9,7}\leq K_{0}\quad \text{for all}\, n\geq 0
\label{9.19}\end{eqnarray}
and this implies $|\widehat{E}|_{9,7}\leq  K_{0}$ because 
each coefficient of $\widehat{E}$ is also that of some $x_{n}$ if $n$
is sufficiently large.

We prove (\ref{9.19}) by induction. For $n=0$ it is true. Assume that
$|x_{n}|_{9,7}\leq  K_{0}$,  for this 
Let us take again the equation (\ref{9.18}) with (\ref{9.2}). We obtain
\begin{eqnarray}
x_{n+1}&+&\widehat{Q}_{1}\widehat{\mathcal{L}}^{-1}\bigg(\frac{\widehat{\mu}}{z^{2}}\bigg)=\widehat{Q}_{1}\widehat{\mathcal{L}}^{-1}\Bigg(\frac{\widehat{W}_{5}}{z^{2}}\widehat{\mathcal{C}}\bigg(\frac{x_{n}}{\widehat{Q}_{1}}\bigg)\Bigg)-\widehat{Q}_{1}\widehat{\mathcal{L}}^{-1}\Bigg(\frac{\widehat{W}_{1}}{2z^{2}\widehat{Q}}\mathcal{Y}_{n}\Bigg)\nonumber\\
&+&\widehat{Q}_{1}\widehat{\mathcal{L}}^{-1}\Bigg(\frac{\widehat{W}_{4}}{2z^{2}\widehat{Q}^{2}}\widehat{\mathcal{C}}\bigg(\frac{x_{n}}{\widehat{Q}_{1}}\bigg)\mathcal{Y}_{n}^{2}\Bigg)
- \widehat{Q}_{1}\widehat{\mathcal{L}}^{-1}\Bigg(\frac{\widehat{W}_{2}}{2z^{2}}\widehat{\mathcal{C}}_{2}\bigg(\frac{\mathcal{Y}_{n}}{\widehat{Q}}\bigg)\Bigg)
\label{9.20}\end{eqnarray}
where $\mathcal{Y}_{n}=\widehat{\mathcal{C}}^{-1}(x_{n})$.
Then using (\ref{9.2})..(\ref{9.17}), we obtain 
\begin{eqnarray*}
\Bigg|\widehat{Q}_{1}\widehat{\mathcal{L}}^{-1}\Bigg(\frac{\widehat{W}_{5}}{z^{2}}\widehat{\mathcal{C}}\bigg(\frac{x_{n}}{\widehat{Q}_{1}}\bigg)\Bigg)\Bigg|_{11,10}&\leq&
K_{1}|x_{0}|_{9,7},\\
\Bigg|\widehat{Q}_{1}\widehat{\mathcal{L}}^{-1}\Bigg(\frac{\widehat{W}_{1}}{2z^{2}\widehat{Q}}\mathcal{Y}_{n}\Bigg)\Bigg|_{17,10}&\leq& 
K_{2}|x_{0}|_{9,7},\\
\Bigg|\widehat{Q}_{1}\widehat{\mathcal{L}}^{-1}\Bigg(\frac{\widehat{W}_{2}}{2z^{2}}\widehat{\mathcal{C}}_{2}\bigg(\frac{\mathcal{Y}_{n}}{\widehat{Q}}\bigg)\Bigg)\Bigg|_{17,10}&\leq&
K_{3}|x_{0}|_{9,7},\\
\Bigg|\widehat{Q}_{1}\widehat{\mathcal{L}}^{-1}\Bigg(\frac{\widehat{W}_{4}}{2z^{2}\widehat{Q}^{2}}\widehat{\mathcal{C}}\bigg(\frac{x_{n}}{\widehat{Q}_{1}}\bigg)\mathcal{Y}_{n}^{2}\Bigg)\Bigg|_{15,10}&\leq& 
K_{4}|x_{0}|_{9,7}^{2}
\end{eqnarray*}
where
\begin{eqnarray*}
K_{1}:&=&
\pi^{2}M_{9,8}\bigg(1+\frac{3}{2}z^{2}\bigg)L_{10,8}M_{7,5}\bigg(\frac{\widehat{W}_{5}}{z^{5}}\bigg)M_{9,7}\bigg(\frac{z^{2}}{\widehat{Q}_{1}}\bigg)\Big\|\cosh\big(\frac{z}{2}\big)\Big\|_{7,5},\\
K_{2}:&=&\pi^{9}M_{15,8}\bigg(1+\frac{3}{2}z^{2}\bigg)L_{16,8}M_{9,1}\bigg(\frac{\widehat{W}_{1}}{2z^{9}\widehat{Q}}\bigg)\,\,c_{9,7},\\
K_{3}:&=&\pi^{9}M_{15,8}\bigg(1+\frac{3}{2}z^{2}\bigg)L_{16,8}M_{7,-1}\bigg(\frac{\widehat{W}_{2}}{2z^{11}}\bigg)\big\|\cosh(z)\big\|_{7,-1}M_{9,1}\bigg(\frac{z^{2}}{\widehat{Q}}\bigg)\,\,c_{9,7},\\
K_{4}:&=&\frac{1}{56\pi^{2}}M_{13,8}\bigg(1+\frac{3}{2}z^{2}\bigg)L_{14,8}M_{14,8}\bigg(\frac{\widehat{W}_{4}}{2z^{2}}\bigg)M_{18,12}\bigg(\frac{z^{4}}{\widehat{Q}^{2}}\bigg)\alpha_{8,8}\,\,\,c_{9,7}^{2}.
\end{eqnarray*}
This imply 
\begin{eqnarray*}
\bigg|x_{n+1}+\widehat{Q}_{1}\widehat{\mathcal{L}}^{-1}\bigg(\frac{\widehat{\mu}}{z^{2}}\bigg)\bigg|_{11,10}\leq 
\big(K_{1}+K_{2}+K_{3}+K_{4}K_{0}\big)K_{0}
\end{eqnarray*}

Evaluation of these quantities yields approximately

$$\begin{array}{|c|c|c|c|} \hline

K_{1}&K_{2}&K_{3} &K_{4}\hl

 375142.8501573\  & 633.1622891589 &10214.44411459
 &0.000011150281\hl

\end{array}
$$
\\
and altogether
\begin{eqnarray}
\bigg|x_{n+1}+\widehat{Q}_{1}\widehat{\mathcal{L}}^{-1}\bigg(\frac{\widehat{\mu}}{z^{2}}\bigg)\bigg|_{11,10}\leq 
118896679183.
\label{9.21}\end{eqnarray}

Now we consider $\text{g}(z):=z^{2}\widehat{\mathcal{L}}^{-1}\bigg(\frac{\widehat{\mu}}{z^{2}}\bigg)$
 , letting  $z^{2}y=z^{-2}\text{g}(z)$, with (\ref{9.15}) we obtain 
\begin{eqnarray}
\widehat{D}y=\widehat{D}\big(z^{-4}\text{g}\big)=f\big(\widehat{D}\big)\bigg[\frac{\widehat{\mu}(z)(4+z^{2})}{4z^{4}\widehat{V}_{1}}\bigg],
\label{9.22}\end{eqnarray}
where
\begin{eqnarray*}
f(z):=\frac{z}{\sinh(z/2)}.
\end{eqnarray*}  

The coefficient of $z^{k}$ in $f(z)$ is smaller than $4(2\pi)^{k}$  we 
obtain by comparison of the coefficients of $z^{j}$ in (\ref{9.22}) 
\begin{eqnarray*}
(j-1)|\text{g}_{j+3}|\leq 4 \sum_{k=4}^{j}*\,(2\pi)^{k-j}\,h_{k}\frac{(j-1)\,!}{(k-1)\,!}
\end{eqnarray*} 
and hence
\begin{eqnarray*}
\frac{|\text{g}_{j}|\pi^{j}2^{j}}{(j-2)\,!}\leq 4
\sum_{k=4}^{\infty}\,\frac{h_{k}(2\pi)^{k}}{(k-1)\,!}\leq 3668333.
\end{eqnarray*} 
We obtain for add $j\geq 9$
\begin{eqnarray*}
\frac{|\text{g}_{j}|}{(j-10)\,!\pi^{-j}}\leq \pi^{3}\cdot3668333\cdot(j-9)(j-8)(j-7)(j-6)(j-5)2^{-j+3}.
\end{eqnarray*}  
and
\begin{eqnarray*}
\frac{|\text{g}_{j-2}|}{(j-10)\,!\pi^{-j}}\leq \pi^{5}\cdot3668333\cdot(j-9)(j-8)(j-7)2^{-j+5}.
\end{eqnarray*}but,\quad
$\Big\{\frac{\widehat{Q}_{1}}{z^{2}}\text{g}(z)\Big\}_{j}=\text{g}_{j}+\frac{3}{2}\text{g}_{j-2}$, 
with (\ref{9.20}) this yields eventually 
\begin{eqnarray*}
\frac{\big|\{x_{n+1}\}_{j}\big|}{(j-10)\,!\pi^{-j}}\leq
118896679184\quad \text{for}\,\,\, j\geq 81.
\end{eqnarray*}
This implies   
\begin{eqnarray*}
\frac{\big|\{x_{n+1}\}_{j}\big|}{(j-7)\,!\pi^{-j}}\leq
\frac{118896679184}{72\cdot73\cdot74}\leq 305691 \quad \text{for}\,\,\, j\geq 81
\end{eqnarray*}
and with $x_{n+1}\equiv x_{n} \mod z^{81}$  we finally proved that
$|x_{n+1}|_{9,7}\leq K_{0}=308027.35$. Thus the proof of (\ref{9.18}) is
complete.
 
\begin{theo}
\begin{eqnarray*}
|\alpha_{n}|\leq \frac{98048.15}{(n-1)(n-2)(n-3)(n-4)(n-5)(n-6)}\quad \text{for}\,\,\, n\geq 9.
\end{eqnarray*}
and even better
\begin{eqnarray*}
|\alpha_{n}|\leq 37845988419 \frac{(n-10)\,!}{(n-1)\,!}\quad
\text{for odd}\,\,\, n\geq 81.
\end{eqnarray*}
\end{theo}
As a consequence we obtain estimates for $\alpha$
\begin{eqnarray*}
\Bigg|\alpha-\frac{1}{\pi}\summ{n}{1}{N}\alpha_{n}\Bigg|&\leq&
12046752267\sum_{n=N}^{\infty}\frac{(n-10)\,!}{(n-1)\,!}\quad
\text{for}\,\,\, n\geq 9,\\
&\leq& 12046752267\cdot
\frac{(N-10)\,!}{16(N-2)\,!}\bigg(1+\frac{8}{N-1}\bigg),\\
&\leq&\frac{762333542}{N(N-1)(N-2)(N-3)(N-4)(N-5)(N-6)(N-7)}
\end{eqnarray*}
We evaluated $\alpha_{n}$ for $ n=9,11,...,81$ and  obtain 
\begin{eqnarray}
|\alpha-1.264150331|\leq 6\cdot 10^{-7}.
\label{9.23}\end{eqnarray}


\begin{thebibliography}{02}
\bibitem{WE}
{W. Eckhaus}, Asymptotic Analysis of Singular perturbations,
North-Holland, Amsterdam (1979).
\bibitem{FS}
{A. Fruchard, R. Sch\"afke}, Exponentially small splitting of
separatrices for difference equations with small step size,  Journal
of Dynamical and Control Systems, Volume 2, Number 2 / April, 1996,
193-238.
\bibitem{HM}
{V. Hakim, K. Mallick}, Exponentially small splitting of separatrices,
matching in the complex plan and Borel summation, Nonlinearity 6(1993) 
57-70.
\bibitem{LS}
{V.F. Lazutkin,I.G. Schachmannski and M.B. Tabanov}, Splitting of
separatrices for standard and semistandard mappings, Physica D,
40 235-248, (1989).

\bibitem{SV}
{R. Sch\"afke, H. Volkmer}, Asymptotic analysis of the equichordal problem,
J.reine u. Angew. Math. 425(1992), 9-60.
\end{thebibliography}
\end{document}